\documentclass[10pt,a4paper,final]{iopart}

\usepackage[breaklinks=true,colorlinks=true,linkcolor=blue,urlcolor=blue,citecolor=blue]{hyperref}

\expandafter\let\csname equation*\endcsname\relax
\expandafter\let\csname endequation*\endcsname\relax

\usepackage{amsmath}
\usepackage{amssymb}
\usepackage{etoolbox}

\patchcmd{\subequations}{\alph{equation}}{\textit{\alph{equation}}}{}{}

\usepackage{amsfonts}
\usepackage{siunitx}
\usepackage[dvipsnames]{xcolor}

\usepackage[capitalize,nameinlink]{cleveref}

\usepackage{changepage}

\Crefname{pluralequation}{Eqs.}{Eqs.}
\Crefformat{pluralequation}{#2Eqs.~(#1)#3}

\newcommand{\taum}{\ensuremath{\tau_\mathrm{m}}}
\newcommand{\taud}{\ensuremath{\tau_\mathrm{d}}}
\newcommand{\tauf}{\ensuremath{\tau_\mathrm{f}}}

\newcommand{\Iext}{\ensuremath{I_1}}

\newcommand{\dd}{\ensuremath{\mathrm{d}}}

\newcommand{\nspikes}{\ensuremath{n_\mathrm{s}}}

\newcommand{\dotsl}[1]{\ensuremath{#1^\prime}}
\newcommand{\dotfa}[1]{\ensuremath{\dot{#1}}}

\newcommand{\tslow}{\ensuremath{\tau}}

\newcommand{\fastrhs}{\ensuremath{\mathbf{F}}}
\newcommand{\slowrhs}{\ensuremath{\mathbf{G}}}

\newcommand{\Xf}{\ensuremath{\mathbf{X}_\mathrm{f}}}
\newcommand{\Xs}{\ensuremath{\mathbf{X}_\mathrm{s}}}

\newcommand{\dfXf}{\ensuremath{\dotfa{\mathbf{X}}_\mathrm{f}}}
\newcommand{\dfXs}{\ensuremath{\dotfa{\mathbf{X}}_\mathrm{s}}}

\newcommand{\dsXf}{\ensuremath{\dotsl{\mathbf{X}}_\mathrm{f}}}
\newcommand{\dsXs}{\ensuremath{\dotsl{\mathbf{X}}_\mathrm{s}}}

\newcommand{\cmf}{\ensuremath{S_0}}

\newcommand{\FSmf}{\ensuremath{\mathcal{M}_{\rm FS}}}
\newcommand{\JOmf}{\ensuremath{\mathcal{M}_{\rm JO}}}

\newcommand{\LCmfa}{\ensuremath{\mathcal{P}^{\rm a}}}
\newcommand{\LCmfr}{\ensuremath{\mathcal{P}^{\rm r}}}
\newcommand{\LCmf}{\ensuremath{\mathcal{P}}}

\newcommand{\rvmf}{\ensuremath{\mathcal{M}_{rv}}}

\newcommand{\FL}{\ensuremath{F_{\rm L}}}
\newcommand{\FU}{\ensuremath{F_{\rm U}}}

\newcommand{\HBL}{\ensuremath{H_{\rm L}}}
\newcommand{\HBU}{\ensuremath{H_{\rm U}}}

\newcommand{\FLc}{\ensuremath{\mathcal{F}_{\rm L}}}
\newcommand{\FUc}{\ensuremath{\mathcal{F}_{\rm U}}}

\newcommand{\HBLc}{\ensuremath{\mathcal{H}_{\rm L}}}
\newcommand{\HBUc}{\ensuremath{\mathcal{H}_{\rm U}}}

\newcommand{\eps}{\varepsilon}

\definecolor{mblue}{HTML}{3fAAD9}

\newcommand{\blue}[1]{\textcolor{blue}{#1}}
\newcommand{\olive}[1]{\textcolor{olive}{#1}}

\usepackage[textwidth=37mm, colorinlistoftodos,prependcaption,textsize=small,backgroundcolor=blue!25]{todonotes}

\usepackage{lipsum}
\usepackage{physics}

\usepackage[numbers, sort&compress]{natbib}
\usepackage{titlesec}
\titleformat{\subsubsection}[runin]
       {\bfseries}{\thesubsubsection}{0.5em}{}

\usepackage{geometry}
\reversemarginpar % put todonotes on the left side

\begin{document}

\title{Bursting in a next generation neural mass model with synaptic dynamics: a slow-fast approach}
\author{Halgurd Taher$^{1}$}
\address{$^{1}$Inria Sophia Antipolis M\'{e}diterran\'{e}e Research Centre, 2004 Route des Lucioles, 06902 Valbonne, France}
\ead{halgurd.taher@inria.fr}

\author{Daniele Avitabile$^{1,2}$}
\address{$^{2}$Department of Mathematics, Faculteit der Exacte Wetenschappen, Vrije Universiteit (VU University Amsterdam), De Boelelaan
	1081a, 1081 HV Amsterdam, Netherlands}
\ead{d.avitabile@vu.nl}

\author{Mathieu Desroches$^{1}$}
\ead{mathieu.desroches@inria.fr}

\vspace{10pt}
\begin{indented}
	\item[]	Date: \today
\end{indented}

\begin{abstract}
	We report a detailed analysis on the emergence of bursting in a recently developed neural mass model that takes short-term synaptic plasticity into account. Neural mass models are capable of mimicking the collective dynamics of large scale neuronal populations in terms of a few macroscopic variables like mean membrane potential and firing rate. The one being used here particularly important, as it represents an exact meanfield limit of synaptically coupled quadratic integrate \& fire neurons, a canonical model for type I excitability. In absence of synaptic dynamics, a periodic external current with a slow frequency $\eps$ can lead to burst-like dynamics. The firing patterns can be  understood using techniques of singular perturbation theory, specifically slow-fast dissection. In the model with synaptic dynamics the separation of timescales leads to a variety of slow-fast phenomena and their role for bursting is rendered inordinately more intricate. Canards are one of the main slow-fast elements on the route to bursting. They describe trajectories evolving nearby otherwise repelling locally invariant sets of the system and are found in the transition region from subthreshold dynamics to bursting. For values of the timescale separation nearby the singular limit $\eps=0$, we report peculiar jump-on canards, which block a continuous transition to bursting. In the biologically more plausible regime of $\eps$ this transition becomes continuous and bursts emerge via consecutive spike-adding transitions. The onset of bursting is of complex nature and involves mixed-type like torus canards, which form the very first spikes of the burst and revolve nearby fast-subsystem repelling limit cycles. We provide numerical evidence for the same mechanisms to be responsible for the emergence of bursting in the quadratic integrate \& fire network with plastic synapses. The main conclusions apply for the network, owing to the exactness of the meanfield limit.
\end{abstract}

%\maketitle

\section{Introduction}\label{sec:intro}
%	In physics and nonlinear dynamics the interaction of particles or nodes within a network can lead from high 	\htnote{This paragraph can be removed probably, it was intended as a link to electrophysiology, local field potentials, EEG,...; but let's see} dimensional microscopic dynamics to low dimensional macroscopic behavior. These problems can typically be treated within the thermodynamic limit, that can lead to a meanfield formulation of the problem under certain assumptions. Meanfield approaches are powerful elements to answer experimental questions, whenever the microscopic scale is inaccessible, while effects on the macroscopic one can be measured.
%	Neuroscience, specifically electrophysiology,  has proven to be in need of such meanfield formulations. 
In the past decade a novel approach in meanfield theory, the so-called Ott-Antonson (OA) ansatz \cite{ottLowDimensionalBehavior2008a,ottLongTimeEvolution2009}, has received major attention. The ansatz serves as a recipe to perform an exact reduction from a high dimensional network of interacting phase oscillators towards a low dimensional dynamical system, that describes the macroscopic behavior. It was firstly applied to the prototypical model for synchronization, namely the Kuramoto model \cite{kuramotoInternationalSymposiumMathematical1975}. The term \textit{exact} refers to the property of the meanfield limit to rely on first principles. Thus, in the thermodynamic limit one obtains an exact agreement of the meanfield system with respect to the macroscopic dynamics of the underlying network.

Remarkably, the OA ansatz is applicable to a broader class of phase oscillator networks and thereby has found its way into the field of computational neuroscience. One system of interest is the Quadratic Integrate \& Fire (QIF) neuron, which represents a canonical model for Saddle-Node on Invariant Circle (SNIC) bifurcations and so-called type I excitability.	Given this canonical role of the QIF neuron, it is standing to reason to investigate the dynamics emerging on macroscopic scale. Here the OA ansatz comes into play. Under certain assumptions the QIF neuron is rendered equivalent to the Ermentrout-Kopell model, also referred to as $\theta$-model \cite{ermentroutParabolicBurstingExcitable1986}. It represents a phase oscillator model for neuronal dynamics, which indeed fits into the class of problems that can be treated within the OA framework.
The OA ansatz applied to an ensemble of QIF neurons, in the work of Montbrió, Pazó and Roxin (MPR), has lead to an upsurge of next generation neural mass models \cite{montbrioMacroscopicDescriptionNetworks2015}.

These models try to capture the macroscopic dynamics of networks of spiking neurons, using just a few variables, like here, the population firing rate and mean membrane potential. They hence contribute crucially
to studies of collective phenomena in these otherwise high dimensional dynamical systems and allow for simple analytical and numerical treatment. Various applications of the MPR model has been studied in recent years. They range from the inclusion of delayed synaptic interactions \cite{pazoQuasiperiodicPartialSynchronization2016}, giving rise to chaos,  to studies of cortical oscillations in multipopulation models \cite{schmidtNetworkMechanismsUnderlying2018b} and cross-frequency coupling \cite{segneriThetanestedGammaOscillations2020a}. While the original MPR model accounts for chemical synapses, the methodology can straightforwardly been applied to include as well electrical synapses, formed by gap junctions between neurons \cite{pietrasExactFiringRate2019,montbrioExactMeanfieldTheory2020}. Extensions of the MPR model towards sparse networks \cite{divoloTransitionAsynchronousOscillatory2018} and fluctuation driven dynamics have been proposed \cite{goldobinReductionMethodologyFluctuation2021c}.

These examples of QIF networks and their neural mass counterparts can give rise to interesting dynamical regimes, typically evoked by bistability in the system. Indeed, the original MPR model exhibits a parameter regime where a stable node and a focus coexist \cite{montbrioMacroscopicDescriptionNetworks2015}. This is particularly relevant for the macroscopic response of neuronal ensembles when they are subject to an external current: bistability implies that a time dependent external drive can lead to interesting firing rhythms.

A recent study takes into account synaptic dynamics in form of exponentially decaying action potentials \cite{avitabileCrossscaleExcitabilityNetworks2021}.
In the present work however, we want to explore a QIF network that accounts for synaptic dynamics in form of short-term synaptic plasticity (STP), thus adding to the biological plausibility.	According to the phenomenological STP model of Tsodyks and Markram  two opposing effects must be distinguished: depression, i.e. weakening, and facilitation, i.e. strengthening, of synaptic connections \cite{tsodyksNeuralCodeNeocortical1997}. To our knowledge, the role of STP in next generation neural mass models has received just little attention, despite being highly relevant in neuroscience.  Previous macroscopic models of STP typically make use of the Wilson-Cowan (WC) model, hence are of heuristic nature \cite{wilsonExcitatoryInhibitoryInteractions1972,tsodyksNeuralNetworksDynamic1998}. Nevertheless they helped to develop a novel synaptic theory of working memory \cite{mongilloSynapticTheoryWorking2008}, a cognitive system for short-term information storage and manipulation in the brain. The synaptic theory of WM has triggered various theoretical studies, which use the WC model with STP in multipopulation topologies in order to understand basic WM operations like information loading and recall, as well as to estimate the maximum WM capacity \cite{miSynapticCorrelatesWorking2017,trubutschekTheoryWorkingMemory2017}.

In a recent study an extension of the MPR firing rate equations towards STP was proposed, in order to model WM \cite{taherExactNeuralMass2020}. The meanfield limit, in presence of STP, remains exact. Therefore one can exploit this limit, in order to get insight into the emergence of firing patterns in the network. An aspect that can lead to complex behavior is the timescale separation, which comes along with STP. Depression and facilitation might indeed act on different timescales. As an example, measurements in the prefrontal cortex suggest that the facilitation of synapses can be maintained for seconds, while depression decays within a few hundred milliseconds \cite{wangHeterogeneityPyramidalNetwork2006}.

Synaptic dynamics and additional timescales enrich the dynamical landscape, by giving rise to bistability involving limit cycles \cite{taherExactNeuralMass2020}. This is the foundation for bursting rhythms to emerge. Bursting refers to dynamics that alternates between a quiescent phase and rapid oscillations. When slowly forcing the population of QIF neurons, by virtues of a slowly drifting external current, the system can transit from a quasi-static motion to rapid oscillations associated with the presence of stable cycles in the system with constant external current.

Bursting has been found in various experimental studies in neuroscience
\cite{
	adamsGenerationModulationEndogenous1985,
	connorsIntrinsicFiringPatterns1990a,
	grayChatteringCellsSuperficial1996,
	schwindtQuantitativeAnalysisFiring1997,
	suExtracellularCalciumModulates2001,
	amirOscillatoryMechanismPrimary2002,
	wellmerLonglastingModificationIntrinsic2002,womackActiveContributionDendrites2002}
and theoretical approaches
\cite{
	plantMechanismUnderlyingBursting1975,
	hindmarshModelNeuronalBursting1984,
	rinzelBurstingOscillationsExcitable1985,
	rinzelFormalClassificationBursting1987,
	rinzelFormalClassificationBursting1987b,
	izhikevichNeuralExcitabilitySpiking2000}
not only aim at classifying the observed dynamics, but also mimicking and revealing the mechanisms responsible for the emergence of bursting. %Ideally, simplified models of complex biological processes can result in predictions and inspire new experiments.
While bursting in spiking neural networks is subject of recent studies \cite{bacakMixedmodeOscillationsPopulation2016,ersozCanardinducedComplexOscillations2020}, the mechanisms responsible for their emergence often remain unclear: exploring the state space of large scale networks is tedious and the addition of slow-fast aspects complicates the problem. The exactness of the MPR model helps to overcome this limitation: analytic tools and bifurcation analysis applied to the neural mass model allow to draw conclusions for the microscopic network.

The main results of this work are related to the emergence of bursting in a QIF network with STP. In particular, we investigate the transition from subthreshold oscillations to bursting in presence of an external slow and periodic current. The forcing introduces a clear timescale separation into the problem, giving rise to intricate slow-fast phenomena and allowing for the application of slow-fast dissection methods, to be described later.
As an outlook, the findings comprise a differentiation of the route to bursting, depending on the timescale separation.
For strongly separated timescales, far away from biologically plausible scenarios, the route is complicated, possibly discontinuous in parameter space and it is related to a certain type of so-called \textit{canards}. However, moderate timescale separation reveals a number of intermingled slow-fast mechanisms that lead to a continuous transition from subthreshold oscillations to bursting and are related to different types of canards. Our results are supported by slow-fast arguments and numerical evidence.

\begin{figure}[h!]
	\centering
	\includegraphics[width=1\linewidth]{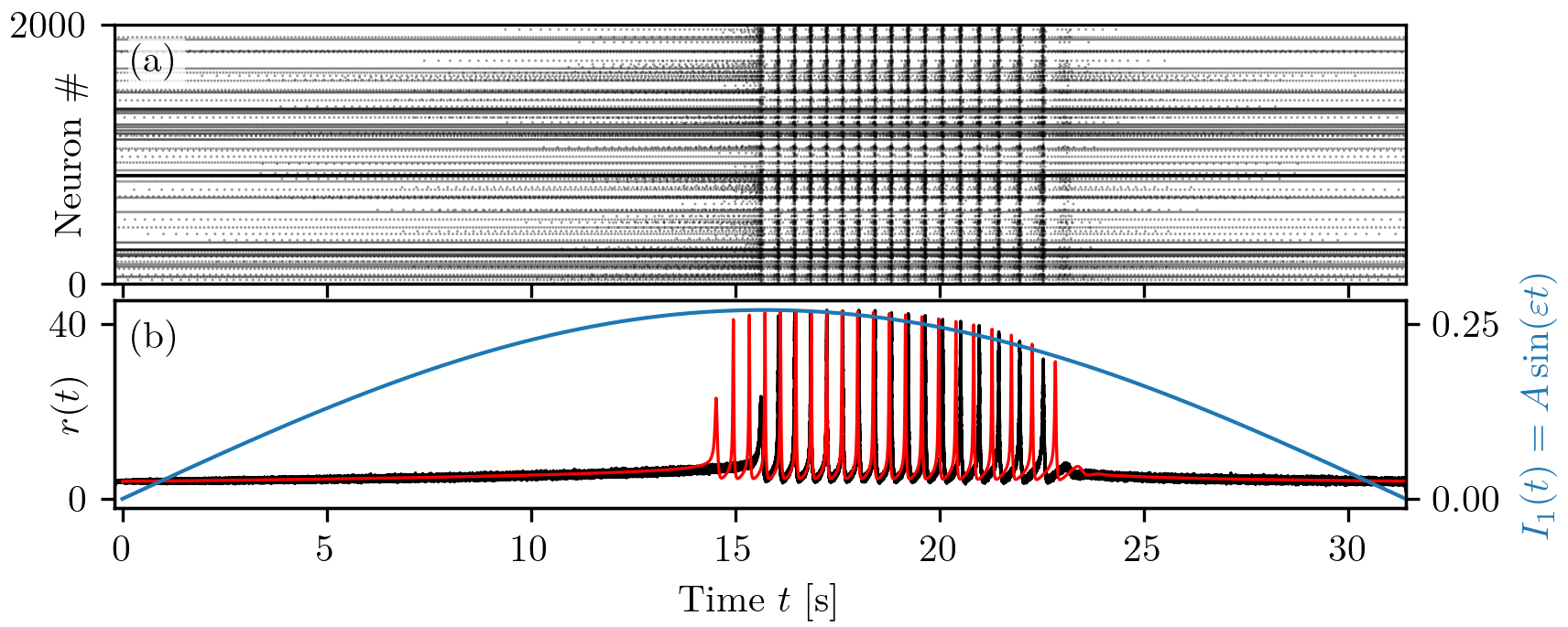}
	\caption{Spiking neuron network and meanfield limit: {\bf (a)} Scatter plot of a subset of 2000 representative neurons out of $10^5$. Each dot represent a spike. {\bf (b)} Firing rate of the network (black) and of the corresponding meanfield limit (red). The blue curve shows a time dependent external current applied to the two systems and is of sinusoidal form.}
	\label{fig:burst_intro}
\end{figure}

A first illustration of the dynamical regime of interest in this work is displayed in \cref{fig:burst_intro}. Panel (a) depicts the response of a large scale network consisting of $N=10^5$ QIF neurons to a slow  external sinusoidal current. The second panel (b) shows the firing rate of the QIF network, as well as the firing rate of the meanfield limit. Both systems undergo a quiescence phase of low firing activity. When the external current exceeds a certain level, the systems start to burst, characterized by a rapid series of synchronized firing at high rates. %The discrepancies between spiking network and meanfield limit arise due to finite size and numerical integration errors,  as well estimations, that are applied in order to simulate the network.

In order to understand how these bursts emerge, we have to encapsulate two main aspects. First, in the upcoming \cref{sec:model},  we will introduce the QIF network model with STP as well as the corresponding meanfield limit, and we will analyze the state space structure and dynamics. Second, the presence of a slow external drive calls for the application of slow-fast dissection. For this it is essential to introduce the general slow-fast framework and revise well known slow-fast mechanisms, which play a role for our model, in \cref{sec:sfframework}. As a next step, this generic methodology for timescale separated problems is applied to the present model in  \cref{sec:dissection}. Dissection is crucial for understanding the results in \cref{sec:fsdynamics}, where canards and in particular jump-on canards are studied using slow-fast arguments. This paves the way to investigate the mechanisms responsible for the emergence of bursting, as done in \cref{sec:transbursting}. Finally, in \cref{sec:netbehaviour} the initially posed problem of bursting on macroscopic scale is approached by a comparison of meanfield dynamics versus QIF network dynamics in the bursting regime.

\section{A next-generation neural mass model with synaptic dynamics}\label{sec:model}

Spiking neuron models can be characterized by their response to the injection of a current, which is often measured in terms of the $f$-$I$ curve, determining the relation of firing frequency $f$ versus input current $I$.
The dynamics of Hodgkin-Huxley type neurons can either be in the excitable or tonic regime \cite{hodgkinLocalElectricChanges1948}. Excitable neurons in absence of input approach an equilibrium. However, sufficient input can excite the membrane potential beyond a threshold leading to the firing of a single action potential, before going back to the rest state. Tonic neurons on the other hand fire periodically with a frequency $f$. Based on the behaviour at the transition from excitable to tonic dynamics one can distinguish (at least) two classes of membranes. For class I neurons the $f$-$I$ curve is continuous and transitions from quiescence ($f=0$) to repetitive firing at arbitrarily slow frequencies ($f>0$). Typically it occurs at a SNIC bifurcation. Class II neurons on the other hand exhibit a discontinuous $f$-$I$ curve, leading to finite firing rates at the onset of tonic firing, and they are usually associated with a Hopf bifurcation. This Hopf bifurcation is often subcritical, for example in the Hodgkin-Huxley, FitzHugh-Nagumo \cite{fitzhughImpulsesPhysiologicalStates1961,nagumoActivePulseTransmission1962} and Morris-Lecar model \cite{morrisVoltageOscillationsBarnacle1981}.

\subsection{Network of spiking neurons and meanfield limit}
The model of interest for our work the canonical model for type I excitability: the QIF neuron. In a network of $N$ synaptically coupled neurons the membrane potentials $V_i(t)$ obey \cref{eq:QIFnetwork1}.
\begin{subequations}\label[pluralequation]{eq:QIFnetwork1}
	\begin{align}
		\dotfa{V}_i & =V_i^2 + \eta_i + J r(t) + I_1(t)                                    \\
		            & \text{if $V_i>V_{\rm thresh}$: $V \leftarrow V_{\rm r}$}\nonumber    \\
		r(t)        & =\frac{1}{N}\sum_{j=1}^{N} \sum_{k: t_j^{(k)}<t} \delta(t-t_j^{(k)})
	\end{align}
\end{subequations}
The total current applied to the neuron is a sum of the constant component $\eta_i$, the synaptic current $Jr(t)$, with synaptic weight $J$ and an external, possibly time-dependent, current $I_1(t)$. Variable $r$ denotes the instantaneous firing rate and is composed of the single neuron spike trains $\sum_{k: t_j^{(k)}<t} \delta(t-t_j^{(k)})$, where $t_j^{(k)}$ denotes the $k$-th spike time of neuron $j$ entering into the Dirac $\delta$ function. Whether a QIF neuron is excitable or tonic depends on $\eta_i$. Given $r=0$,  the neuron is excitable for $\eta_i<0$ and in tonic firing state for $\eta_i>0$.
Firing occurs whenever $V_i$ exceeds the threshold $V_{\rm thresh}$ at which the reset rule applies leading to a reset of the potential to $V_{\rm r}$.

When performing the thermodynamic limit $N\to \infty$ and imposing certain conditions, the mean dynamics of the above microscopic model leads to a reduced macroscopic description in terms of the mean membrane potential $v(t)$ and firing rate $r(t)$, namely the MPR model \cite{montbrioMacroscopicDescriptionNetworks2015}.	The derivation is based on the OA ansatz and yields an exact reduction \cite{ottLowDimensionalBehavior2008a}. Thus the collective behavior of the QIF network, aside from finite size fluctuations, is in perfect agreement with the MPR model.
Following the OA Ansatz and the derivation of the MPR model the following assumptions have to be made in order to obtain an exact firing rate formulation.

(i) The threshold and reset voltage have to be considered in the limit $V_{\rm thresh}=-V_{\rm r}\rightarrow \infty$, rendering the QIF neuron identical to the $\theta$-model \cite{ermentroutParabolicBurstingExcitable1986}.

(ii) The excitabilities $\eta_i$ are drawn from a Lorentzian distribution $g(\eta)=\frac{1}{\pi}\frac{\Delta}{(\eta-\bar{\eta})^2+\Delta^2}$, centred at $\bar{\eta}$ and with the width parameter $\Delta$.

(iii) Neurons are all-to-all coupled. This way each QIF neuron receives identical synaptic current $Jr(t)$.

(iv) The QIF network has to be considered in the thermodynamic limit $N\rightarrow \infty$.

\noindent The resulting MPR model consists of two ordinary differential equations for $r(t)$ and $v(t)$ given in \cref{eq:montbrioFRE}.
\begin{subequations}\label[pluralequation]{eq:montbrioFRE}
	\begin{align}
		\dotfa{r} & =\frac{\Delta}{\pi} + 2rv                        \\
		\dotfa{v} & =v^2  - (\pi  r)^2 + J r+ \bar{\eta}  + \Iext(t)
	\end{align}
\end{subequations}

Despite having a rather simple state-space structure, the MPR model can give rise to interesting periodic patterns when externally forced. In the case of constant $\Iext$ periodic solutions are absent and one can find node, focus and saddle equilibria. However, there are regions of parameter space in which bistability between the node and focus appears. Hence slow periodic forcing, for example given by $\Iext=A\sin(\varepsilon t), \quad 0<\varepsilon\ll 1$, can lead to a hysteretic loop in these regions, as shown in \cite{montbrioMacroscopicDescriptionNetworks2015}. In this case the trajectories consist of a low firing rate segment and a high firing rate segment with damped oscillations, related to the presence of foci in the system with constant $I_1$. Orbits of this type can already be seen as bursting patterns characterized by an alternation of slow drifts and fast oscillations. However, in the limit $\varepsilon=0$ of infinitely slow forcing, the fast oscillations vanish. In that case, the resulting cycles can be classified as relaxation oscillations, which are introduced in \cref{sec:sfframework}.

\subsection{Neural mass model with short-term plasticity}\label{sec:neurmassSTP}
Our present study investigates bursting patterns in an extended version of \cref{eq:montbrioFRE} that accounts for STP as described in the phenomenological model of Tsodyks and Markram \cite{tsodyksNeuralCodeNeocortical1997}. Short-term synaptic depression is related to neurotransmitter depletion. Each neuron $i$ has a limited amount $X_i(t)\in [0,1]$ of resources (i.e vesicles ready to be released). Spiking is followed by the emittance of presynaptic action potentials. Upon their arrival at the synaptic terminal a fraction $U_i(t)\in [U_0,1]$ of the neurotransmitters is released into the synaptic cleft, resulting in the generation of postsynaptic potentials (PSPs). Therefore, each presynaptic spike is linked to the utilization and reduction of resources available for the generation of upcoming PSPs, consequently leading to a decrease of future postsynaptic excitations. The resource $X_i$ exponentially recovers to its base value of $X_i=1$ on a timescale $\taud=\SI{200}{ms}$ (depression timescale).

Facilitation, as opposed to depression, leads to enhanced PSPs and is related the neurotransmitter release probability at the synaptic terminals, which is modelled by the utilization factor $U_i$.  The release probability (and therefore $U_i$) depends on the intracellular calcium concentration. The neurotransmitter release is associated with the accumulation of calcium ions in the presynaptic terminal, hence each spike leads to an increase of $U_i$. Calcium concentration and the utilization factor decay to the base level $U_i=U_0$ on the facilitation timescale $\tauf=\SI{1500}{ms}$.

We will focus on an implementation of STP into the model on macroscopic level (\textit{m-STP}) as suggested in \cite{taherExactNeuralMass2020}, in order maintain the exactness of the firing rate model. In other words, depression and facilitation, accounted for by $X_i$ and $U_i$ respectively, will not be treated on single neurons level, but rather on population level, with the depression and facilitation variables $x(t)$ and $u(t)$, respectively. This results in $N$ membrane potential equations and two synaptic equations for the QIF network, as given in \cref{eq:QIFnetwork2}.
\begin{subequations}\label[pluralequation]{eq:QIFnetwork2}
	\begin{align}
		\dotfa{V}_i & =V_i^2 + \eta_i + Jux r +\Iext(t)                         \\
		\dotfa{x}   & =	\frac{1-x}{\taud}-uxr\label{eq:QIFdepression}            \\
		\dotfa{u}   & =	\frac{U_0-u}{\tauf}+U_0(1-u)r \label{eq:QIFfacilitation}
	\end{align}
\end{subequations}

The amount of resources $x$ of the QIF network reduces when the population firing rate $r$ increases, while at the same time the utilization $u$ increases. Both quantities enter into the effective synaptic weight $Jux$. %\citeauthor{schmutz2020mesoscopic} classified this STP meanfield limit as a first order limit, which as opposed to the second order limit, neglects the correlation of single neuron 
The extended system given in \cref{eq:model1}, in the following referred to as \textit{neural mass with STP} (NMSTP), represents an exact meanfield limit of the QIF network given in \cref{eq:QIFnetwork2}. The state variables are the firing rate $r$, mean membrane potential $v$,  amount of resources $x$ and utilization factor $u$.
\begin{subequations}\label[pluralequation]{eq:model1}
	\begin{align}
		\dotfa{r} & =\frac{\Delta}{\pi } + 2rv                                         \\
		\dotfa{v} & =v^2  - (\pi r)^2 + Ju x r  + \bar{\eta} + \Iext(t)\label{eq:FREv} \\
		\dotfa{x} & =	\frac{1-x}{\taud}-uxr\label{eq:FREdepression}                     \\
		\dotfa{u} & =	\frac{U_0-u}{\tauf}+U_0(1-u)r \label{eq:FREfacilitation}
	\end{align}
\end{subequations}
We note that \cref{eq:model1} will evolve on the fastest timescale of our problem. This holds despite the fact that it already possesses multiple timescale via $\taud$ and $\tauf$. However, as we will  discuss later, this inherent timescale separation of the NMSTP is subtle and not observable everywhere in state space. Nevertheless, it has significant impact on how the transition from subthreshold (non-bursting) behaviour to bursting occurs (see also \cref{sec:jump-on,sec:blocking}).

\subsection{Dynamics under constant forcing}
\begin{figure}[h]
	\centering
	\includegraphics[width=1\linewidth]{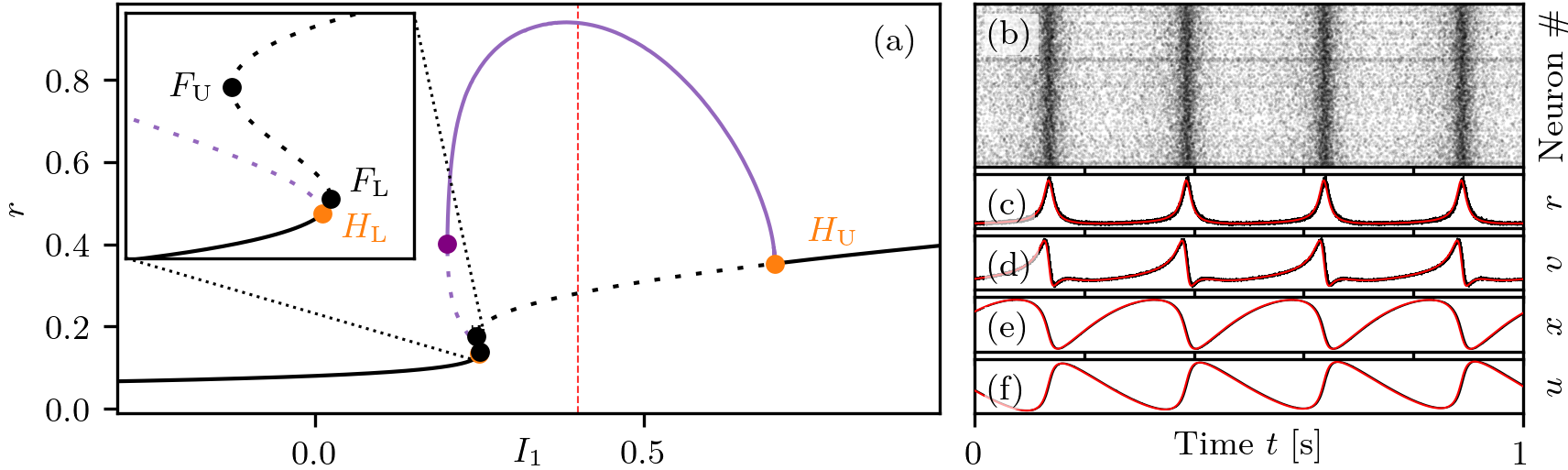}
	\caption{Solution families of the system with constant forcing: {\bf (a)} Bifurcation diagram $r$ versus $I_1$ of \cref{eq:model1}.  For $I_1\lesssim 0.25$ there exists only one fixed point (FP, solid black line). At  $I_1\approx 0.25$ the FP destabilizes via a subcritical Hopf-Bifurcation ($\HBL$, lower orange dot), creating a branch of unstable limit cycles (LC, purple dashed line). Two saddle-node bifurcations (black dots) of the unstable FP branch (dashed black line) occur in a narrow regime of $I_1$, folding the branch twice. Stability is regained for $I_1\gtrsim 0.7$ at a supercritical Hopf-Bifurcation ($\HBU$, upper orange dot). The unstable LC stabilizes (solid purple line) via a saddle-node bifurcation of cycles (purple dot) and vanishes at the second Hopf Bifurcation. The purple line marks the maximum firing rates of the LC branch. {\bf (c - e)} Periodic solution $(r(t),v(t),x(t),u(t))$ vs. time $t$ for $I_1=0.4$ marked in panel (a) by a dashed red line. The red curves show simulations results of the NMSTP, the black ones of the network \cref{eq:QIFnetwork2}. {\bf (b)} Spike scatter plot for 20000 representative neurons in the network out of $N=100000$.}
	\label{fig:bifdiag}
\end{figure}

Most of the parameters values used for \cref{eq:model1} will remain fixed and, if not stated differently, given in \cref{tab:params1}. Note that time is measured in units of the membrane time constant $\taum$. For more details on the numerical methods, we refer to the supplementary material of this work.
\begin{table}[ht]
	\caption{Parameters and their values, which are fixed throughout this work, if not stated differently.}\label{tab:params1}
	\centering
	%		\footnotesize
	\begin{tabular}{@{}lll}
		\br
		Symbol       & Description            & Value                 \\
		\mr
		$\Delta$     & Width of Lorentzian    & 0.5                   \\
		$\bar{\eta}$ & Centre of Lorentzian   & -1.7                  \\
		$J$          & Synaptic weight        & 30                    \\
		$U_0$        & Baseline utilization   & 0.1                   \\
		$\taum$      & Membrane time constant & $\SI{20}{ms}$         \\
		$\taud$      & Depression timescale   & $\SI{200}{ms}/\taum$   \\
		$\tauf$      & Facilitation timescale & $\SI{1500}{ms}/\taum$ \\
		\br
	\end{tabular}
\end{table}
We will outline the different dynamical regimes in presence of a constant current $I_1(t)=\text{const.}$, using the above parameter values.	In \cref{fig:bifdiag}(a) the resulting bifurcation diagram is displayed.

For currents $I_1\lesssim 0.25$ we find a branch of stable node equilibria at low firing rates. The branch develops into a family of foci and destabilizes around $I_1\approx 0.25$ via a subcritical Hopf bifurcation ($\HBL$) followed by two saddle-node (fold) bifurcation at $\FL$ and $\FU$ (black dots), where $F_k=(r_k,v_k, x_k,u_k,I_k),\quad k\in \{\mathrm{L},\mathrm{U}\}$, denotes the equilibrium and parameter values of the bifurcations.
These folds can also be found in absence of STP, in which case the upper branch is stable. However, in \cref{fig:bifdiag}(a) the instability persists throughout the S-shaped curve up until the upper supercritical Hopf bifurcation $\HBU$.
The lower Hopf bifurcation $\HBL$ generates a family of unstable limit cycles that undergoes a fold bifurcation of cycles, giving rise to stable periodic solutions.

One of these trajectories $(r(t),v(t),x(t),u(t))$ is presented in \cref{fig:bifdiag}(b-f) as a function of time. It is superimposed onto the corresponding variables calculated by simulating a QIF network governed by \cref{eq:QIFnetwork2} and consisting of $N=100000$ neurons. For this network, the firing rate is estimated via binning of time, i.e, by counting the number of spikes per time bin of width $\Delta t=10^{-2}$. The average membrane potential for the network reads $v(t)=\frac{1}{N}\sum_{j=1}^N V_j(t)$.

The primary mechanism driving the oscillations is an interplay of so-called populations bursts  and the ensuing synaptic depression and facilitation. At the microscopic scale, population bursts are emitted via a cascade of spikes throughout the network, which as a consequence leads to the facilitation of synapses, leveraging the firing activity further; see \cref{fig:bifdiag}(b,c,f). The consequent depression suppresses the activity, but recovers on the timescale $\taud$ allowing for the emittance of population bursts in a periodic manner.

Notably, in the $I_1$-interval depicted in the inset of \cref{fig:bifdiag}(a) we find bistability between equilibria and limit cycles. We can therefore predict that a time dependent slow current $I_1(t)$, evolving across this region, will lead to a dynamic transition from the equilibrium branch to the stable limit cycles, giving rise to bursting. This exact example can be found in \cref{fig:burst_intro}.

Overall, in contrast to the QIF network without STP and original MPR model, where no limit cycles exist, STP gives rise to bistability among equilibria and cycles. In Ref. \cite{montbrioMacroscopicDescriptionNetworks2015} a slow periodic currents leads to the emergence of macroscopic relaxation-type oscillations in the network. We want to investigate how the presence of STP impacts the response of the system towards such an input. Simulations of QIF networks are difficult computationally. However, the expected agreement of QIF network results and the NMSTP depicted in \cref{fig:bifdiag}(b-f) justifies to perform the upcoming analysis using solely the NMSTP. We will return to the implication of NMSTP dynamics for the network in \cref{sec:netbehaviour}.

\subsection{Dynamics under slow periodic forcing}
Owing to the previous observations in the system with constant forcing, we will introduce a slow periodic drive into the model via the external current $I_1$. We impose that it evolves periodically and on a timescale considerably larger than the slowest timescale of the neural mass, namely the facilitation decay time $\tauf$. %Scenarios in which neuronal populations are periodically driven are ubiquitous. 
In order to remain in a general framework, $\Iext(t)$ will be sinusoidal, given by $\Iext(t)=A\sin{(\eps t)}$, 
with period $T=\frac{2\pi}{\eps}\gg \tauf$ and amplitude $A$. Throughout this work we set $\tauf=\SI{1500}{ms}/\taum=75$, therefore the separation between forcing and slowest intrinsic timescale of the fast subsystem is calculated as $\tauf/T=\eps\frac{\tauf}{2\pi}\approx 10\eps$.
%	\htnote[inline]{examples: cortical oscillations, interactions with other populations, sensory input, homeostatic changes...}
%	\htnote[inline]{Some values: $\eps \in [2\cdot 10^{-4}, 2\cdot 10^{-3}]$ approximately corresponds to $ T \in [ \SI{600}{s},\SI{60}{s}]$. In this range we find canards, jump-on canards, bursting, mixed-type torus canards, basically all phenomena described in this work}
%	\mdnote[inline]{Good :)}

Through the choice of $\Iext$ to be explicitly time dependent, the system given in \cref{eq:model1} becomes non-autonomous. This in turn comes along with hurdles in the application of slow-fast dissection. Thus, in order to retrieve an autonomous system, a second forcing variable $I_2$ is introduced. The dynamics of ($I_1,I_2)$ follows a Hopf normal form as given below.
\begin{subequations}\label[pluralequation]{eq:forcing1}
	\begin{align}
		\dotfa{I_1} & =\eps g_1(I_1,I_2)=\eps\left[I_1(a-I_1^2-I_2^2)+I_2\right] \\
		\dotfa{I_2} & =\eps g_2(I_1,I_2)=\eps\left[I_2(a-I_1^2-I_2^2)-I_1\right]
	\end{align}
\end{subequations}
The Hopf bifurcation at $a=0$ gives rise to stable limit cycles of the form $(I_1,I_2)=A\cdot(\sin {\eps t},\cos{\eps t})$, with amplitude $A=\sqrt{a}$ and angular frequency $\eps$, in the following referred to as \textit{forcing cycle}. To assure equivalence of the explicitly defined $I_1(t)=A\sin (\eps t)$ and the one generated by the Hopf form \cref{eq:forcing1}, the initial conditions $\left((I_1(t_0), I_2(t_0)\right)$ will lie on $(I_1,I_2)=A\cdot(\sin {\eps t},\cos{\eps t})$.
The full system is given by the NMSTP in presence of slow external forcing, i.e, \crefrange{eq:model1}{eq:forcing1}.

\begin{figure}[h]
	\centering
	\includegraphics[width=1\linewidth]{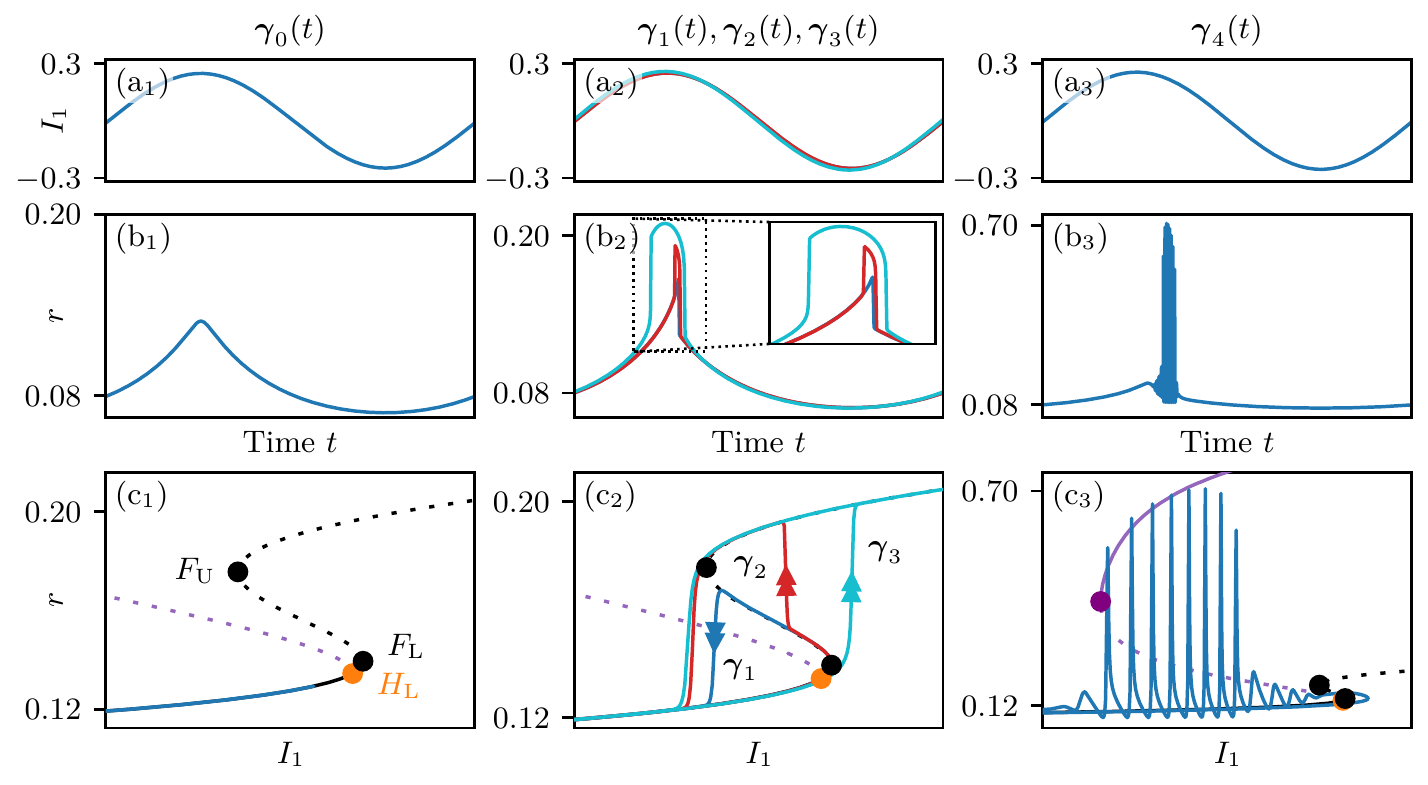}
	\caption{Typical solutions $\boldsymbol{\gamma}_0(t)$ to $\boldsymbol{\gamma}_4(t)$ of the full system: {\bf (a)} Periodic forcing current $I_1(t)$ and {\bf (b)} firing rate $r(t)$ vs. time $t$. {\bf (c)} Same trajectories superimposed on the bifurcation diagram of the unforced system in $r$-$I_1$.
	The parameter values are as follows:
	\mbox{$\eps=10^{-5}$};
	$\boldsymbol{\gamma}_0$:  \mbox{$A\approx 0.2487$};
	$\boldsymbol{\gamma}_1$ to $\boldsymbol{\gamma}_3$: in increasing order exponentially close to \mbox{$A\approx 0.2507$};
	$\boldsymbol{\gamma}_4$: \mbox{$\eps=10^{-3}$}, \mbox{$A\approx 0.2553$}.}
	\label{fig:fullsystemsolutions1}
\end{figure}

To understand the impact of this slow forcing, it is advantageous to superimpose solutions of the full problem on the $r$ vs. $I_1$ bifurcation diagram of the unforced system, this is at the core of the slow-fast dissection introduced by J. Rinzel~\cite{rinzelBurstingOscillationsExcitable1985,rinzelFormalClassificationBursting1987,rinzelFormalClassificationBursting1987b}.
An example of a purely slow trajectory is shown in \cref{fig:fullsystemsolutions1}(a$_1$-c$_1$) and labeled $\boldsymbol{\gamma}_0(t)$. The firing rate $r(t)$, shown in panel (b$_1$), increases and decreases following the same pattern as the forcing $I_1(t)$ in panel (a$_1$). Moreover, in $r-I_1$ projection, shown in panel (c$_1$), it becomes clear that the dynamics takes place nearby the equilibrium branch of the unforced system. The forcing introduces a drift of the equilibrium, slow enough to be followed by the dynamics in an $\mathcal{O(\eps)}$ neighbourhood of the branch.

While this example can be understand as a quasi-static motion, the more complex solutions $\boldsymbol{\gamma}_1(t)$ to $\boldsymbol{\gamma}_4(t)$ in columns 2 and 3 of \cref{fig:fullsystemsolutions1}, exhibit canard dynamics and bursting, respectively. A more rigorous analysis is required, including a slow-fast dissection of the model. In the next sections a generic framework in which timescale separated problems can be treated is introduced and applied to the NMSTP.

%	We will provide a brief overview of solutions of the full system before getting into the details on of how these arise from slow-fast properties.

\section{Slow-fast framework and state of the art}\label{sec:sfframework}	The dynamics of slow-fast systems can be regarded in terms of fast variables  $\Xf(t)\in \mathbb{R}^k$ and slow variables $\Xs(t)\in \mathbb{R}^l$. Their dynamics is governed by the differential equations given in \cref{eq:fullfast} and here referred to as \textit{full system},
\begin{subequations}\label[pluralequation]{eq:fullfast}
	\begin{align}
		\dfXf & =~\fastrhs(\Xf,\Xs)     \\
		\dfXs & =\eps \slowrhs(\Xf,\Xs)
	\end{align}
\end{subequations}
with fast-time parametrisation $t$ (the overdot denoting differentiation with respect to $t$), $\fastrhs(\Xf,\Xs):\mathbb{R}^k \times \mathbb{R}^l \;\to\; \mathbb{R}^k$ and
$\slowrhs(\Xf,\Xs): \mathbb{R}^k \times \mathbb{R}^l \;\to\; \mathbb{R}^l$. Here the separation of timescales is reflected by a small parameter $0<\eps \ll 1$. We will refer to this type of system as $k$-fast $l$-slow system.

A different formulation of the full system is obtained in \cref{eq:fullslow} by parametrizing it in slow time $\tau:=\eps t$.
\begin{subequations}\label[pluralequation]{eq:fullslow}
	\begin{align}
		\eps \dsXf & =\fastrhs(\Xf,\Xs)        \\
		\dsXs      & =\slowrhs(\Xf, \Xs) \quad.
	\end{align}
\end{subequations}
The derivative with respect to slow time $\tslow$ is denoted $\dotsl{(\_)} := \dd/\dd\tslow (\_) = \frac{1}{\eps} \dotfa{(\_)}$. The two representations \cref{eq:fullfast} and \cref{eq:fullslow} are equivalent, but they allow to exploit the premise of slow-fast systems, namely the timescale separation given by a small value of $\eps$, in different ways. It is natural to consider the \textit{singular limit}  $\eps=0$ and take the fast and slow time parametrizations into account. One obtains two different subsystems, which represent a dissection of slow and fast dynamics of the full system.

In the first case we obtain the \textit{fast subsystem} \cref{eq:fastsub}. This limit can be used to understand dynamics of the full system for which $\Xf$ evolves fast and results in the following $k+l$ ODEs, $l$ of which being trivial:
\begin{subequations}\label[pluralequation]{eq:fastsub}
	\begin{align}
		\dfXf & =\fastrhs(\Xf,\Xs)\label{eq:fastsubeq} \\
		\dfXs & =\mathbf{0}
	\end{align}
\end{subequations}
Indeed, the dynamics of the slow variables $\Xs$ is trivial and their value does not change in time. As a matter of fact they can be treated as parameters entering into the dynamics of $\Xf$.

The second limit $\eps\rightarrow 0$, now done in the slow-time parametrization \cref{eq:fullslow}, yields the \textit{slow subsystem}, namely:
\begin{subequations}\label[pluralequation]{eq:slowsub}
	\begin{align}
		\mathbf{0} & =\fastrhs(\Xf,\Xs)\label{eq:slowsubcritmf}        \\
		\dsXs      & =\slowrhs(\Xf,\Xs) \quad. \label{eq:slowsubI1I2}
	\end{align}
\end{subequations}
\cref{eq:slowsub} is also referred to as the \textit{reduced system} and it is represented by a differential-algebraic system, in which the dynamics of the slow variables remains unchanged with respect to the full system and is governed by $\dsXs=\slowrhs(\Xf,\Xs)$. The dynamics of the fast variables on the other hand are hidden within the $k$ algebraic constraints \cref{eq:slowsubcritmf}.
They define the \textit{critical manifold}:
\begin{align}
	\cmf=\{(\Xf,\Xs)~|~\fastrhs(\Xf,\Xs)=0\}\quad ,
\end{align}
usually a $l$-dimensional manifold embedded in $\mathbb{R}^{(k+l)}$.

In the slow subsystem the dynamics of the fast variables is slaved to the slow variables, their relation is given by the critical manifold's equations, which defines the state space of this limiting problem: motion of the slow subsystem takes place on $\cmf$. At the same time points of $\cmf$ correspond to equilibria of the fast subsystem, as is clear from equation~\cref{eq:fastsubeq}. By joining solutions of the different subsystems at specific points, \textit{singular orbits} can be constructed. They are trajectories resulting from the concatenation of slow and fast segments, for which the dynamics is determined by the respective subsystem.

Solutions of systems of the type of~\cref{eq:fullfast} with both slow and fast segments are $\eps$-perturbations of singular orbits and the prototypical slow-fast cycles are \textit{relaxation oscillations}~\cite{vanderpolLXXXVIIIRelaxationoscillations1926}. The way these cycles emerge in parameter space is rather peculiar and involve the famous canard solutions, which we review in the context of the classical van der Pol (VdP) system below.

\subsection*{Classical canards in the van der Pol oscillator.}\label{sec:vdpcanard}
Here we briefly present the essentials of classical canards in the prototypical VdP system given in \cref{eq:vdp}, consisting of one fast variable $x$ and one slow variable $y$.	Once written in first-order form, the system's equations read:
\begin{subequations}\label[pluralequation]{eq:vdp}
	\begin{align}
		\dotfa{x} & =y-\frac{x^3}{3}+x \\
		\dotfa{y} & =\eps(a-x)\quad .
	\end{align}
\end{subequations}
The critical manifold $\cmf$ is one dimensional and S-shaped, parametrized by $x$ as $y (x) =\frac{x^3}{3}-x$. In VdP the critical manifold $\{(x,y);\;y=y(x)\}$ has a local maximum and a local minimum, at $x=\pm 1$, respectively. Therefore $\cmf$ folds twice and has three branches, of which the middle one is repelling and the other two attracting. They correspond to unstable and stable equilibria of the fast subsystem, respectively. The slow nullcline, determined by $\dotfa{y}=0$, is the straight line $\{x=a\}$ and, at $a=\pm 1$, it intersects $S_0$ at the fold points perpendicularly.

Overall, the VdP system, for $0<\eps \ll 1$ and with $a$ as a bifurcation parameter, displays four regimes, from subthreshold oscillations to relaxation oscillations. The former is characterized by a quasi-static motion that takes place entirely nearby one of the stable branches of $\cmf$ with a slow frequency and low amplitude. Relaxation oscillations on the other hand comprise a quasi-static motion nearby both attracting branches of $\cmf$, joined via fast transitions: when the branch becomes repelling at one of the fold points, also called \textit{jump-off} point in this context, a rapid \textit{jump} to the opposed attracting branch of $\cmf$ occurs.

The subthreshold regime terminates at a supercritical Hopf bifurcation, at which stable small Hopf cycles emanate, which are not yet of relaxation type and only exist in an $O(\eps)$ distance from the Hopf point. This is followed by an exponentially narrow parameter interval, for which the orbits grow in an explosive manner when the parameter is varied. This phenomenon is known as a \textit{canard explosion}~\cite{bronsBifurcationsInstabilitiesGreitzer1988} and the associated \textit{canards} separate the small Hopf cycles from relaxation oscillations \cite{benoitChasseAuCanard1981}. They evolve for some time near the repelling sheet of $S_0$, before jumping to one of the attracting sheets. They represent a mechanism which allows to connect the attracting sheet of the critical manifold to the repelling one. In our work, as in many problems involving bursting solutions, canard dynamics play an essential role. They separate the parameter regime for subthreshold oscillations from the one where bursting can occur.

The classical canard can be understood by taking into account the \textit{slow flow}, i.e, the flow of \cref{eq:slowsub}, which describes the slow dynamics on $\cmf$, and will be derived later.
For VdP it reads $(\dotsl{x},\dotsl{y})=(\frac{a-x}{x^2-1},a-x)$ and is in general undefined at the folds ($x=\pm 1$). 	However, its $x$-component reduces to $\dotsl{x}=-\frac{1}{x \pm 1}$ for $a=\pm 1$. In this case the flow remains defined at the corresponding fold $x=\pm a$, but undefined at the other fold $x=\mp a$. Despite the intersection of the slow nullcline $x=a$ with the cubic nullcline, no equilibrium exists. Instead a \textit{turning point} forms, which allows a continuous passage through the corresponding fold without any obstruction. This passage from the attracting to the repelling sheet of $\cmf$ occurs in finite time and is the basis of singular canard orbits. Therefore the turning point is also referred to as \textit{canard point}.

\subsubsection*{Folded-saddle canards.}\label{sec:fscanards}
The understanding of canards based on the VdP system can be extended towards 1-fast 2-slow systems. In the most simple case the parameter $a$ of the VdP system is subject to a slow drift given by $\dot{a}=\eps \mu$, with a constant speed $\mu$. Naturally, the slow nullcline is absent in this case and the intersection with the fast subsystem's fold, which is forming the \textit{turning point}, can not occur. Nevertheless canards can be found and are a result of a slow passage effect: the parameter $a$ dynamically transitions through the canard explosion of the original problem.

In VdP the existence of turning points is conditioned by the fact that the slow variable receives feedback from the fast variable , leading to a ``turn" of the slow flow direction, depending on the value of $x$. In the neural mass model with STP in presence of external periodic forcing, however, this feedback is absent. The turn of the slow flow appears naturally via the form of the external forcing. In both cases the passage through the fold underlies the same mechanisms and it is well understood by making use of an auxiliary system, called the \textit{desingularized reduced system} (DRS), which is introduced in introduced in Section~\nameref{sec:desingularization}

Making use of the DRS, the formation of canards as described above, is reflected by so-called \textit{folded-saddle singularities}, which allow a passage from the attracting to the repelling sheet of $\cmf$.  Ultimately, this leads to the existence of canard orbits in the extended VdP system, as well as in the neural mass with STP in presence of periodic forcing.
%	Making use of the DRS, the formation of canards as described above, is reflected by so\olive{-}called \textit{folded\olive{-saddle} singularities} and in particular the existence of a \textit{folded homoclinic} connection. The latter allows a passage from the attracting to the repelling sheet of $\cmf$.  Ultimately, this leads to the existence of canard orbits in the extended VdP system, as well as in the neural mass with STP in presence of periodic forcing.

\subsubsection*{Torus and mixed-type canards.}
The term canard	is not restricted to dynamics taking place in the vicinity of (or on) attracting and repelling manifolds, which represent equilibria. In general, it refers to any type of solution evolving near attracting and repelling invariant sets associated with the fast subsystem. These invariant sets can correspond to equilibria but also to limit cycles. Following this definition, a particular type of canard can be found in elliptic bursters \cite{baspinarCanonicalModelsTorus2021}, which require at least a 2-fast 1-slow system. Here elliptic bursting can arise due to a subcritical Hopf-Bifurcation (in the fast subsystem) giving rise to unstable limit cycles, which stabilize via a fold bifurcation of cycles. The Hopf bifurcation initiates the burst, while the fold of cycles marks their termination. Usually in elliptic bursters the full dynamics follow the family of stable limit cycles of the fast subsystem. However so\olive{-}called \textit{torus canards} can be found for small enough $\eps$ \cite{benesElementaryModelTorus2011,burkeShowcaseTorusCanards2012}. They describe orbits following a stable family of fast subsystem cycles and switching to the unstable one past the fold.

A hybrid of classical canards and torus canards, so-called mixed-type canards, were reported in \cite{desrochesCanardsMixedType2012}. They describe trajectories that spend time near repelling branches of equilibria as well as limit cycles, and can therefore be seen as a mix of classical canards and torus canards. Segments of these solutions evolve nearby unstable equilibria and connect to unstable limit cycles of the fast subsystem.
These last types of canards are found in bursting systems, which are ubiquitous in the modeling of neural activity at both single-cell and population level.

\section{Slow-fast dissection of the model}\label{sec:dissection}
We will start a systematic investigation of the full system by dissecting it into a slow and fast subsystem. The full problem represents a 4-fast 2-slow system with $\Xf=(r,v,x,u)$ and $\Xs=(I_1,I_2)$. Their dynamics is governed by the right hand sides $\fastrhs(\Xf,\Xs)$ and $\slowrhs(\Xs)$, recalled below.
\begin{align}
	\fastrhs(\Xf,\Xs) & =
	\begin{pmatrix}
		\frac{\Delta}{\taum \pi}+2rv             \\
		v^2  + J\taum uxr - (\pi \taum r)^2 +I_1 \\
		(1-x)/\taud-uxr                          \\
		(U_0-u)/\tauf+U_0(1-u)r
	\end{pmatrix}\label{eq:model2} \\
	\slowrhs(\Xs)     & =
	\eps\begin{pmatrix}
		g_1(I_1,I_2) \\
		g_2(I_1,I_2)
	\end{pmatrix}
	=\eps
	\begin{pmatrix}
		\left(I_1(a-I_1^2-I_2^2)+I_2\right) \\
		\left(I_2(a-I_1^2-I_2^2)-I_1\right)
	\end{pmatrix}\label{eq:forcing2}
\end{align}
The equilibrium branches of the fast subsystem, shown in \cref{fig:bifdiag}(a), are defined via $\{\Xf|\dotfa{\Xf}=0\}$, with the slow variable coordinates acting as bifurcation parameters. This naturally coincides with the definition \cref{eq:critmf} of the critical manifold $S_0$
\begin{align}\label{eq:critmf}
	\fastrhs(\Xf,\Xs) & =\mathbf{0}
\end{align}

%	\begin{subequations}\label{eq:critmf}
%		\begin{align}
%			0&=\frac{\Delta}{\taum \pi}+2rv\\
%			0&=		v^2  + J\taum uxr - (\pi \taum r)^2 +I_1\\
%			0&=	(1-x)/\taud-uxr\\
%			0&=	(U_0-u)/\tauf+U_0(1-u)r
%		\end{align}
%	\end{subequations}

%	\begin{subequations}\label{eq:critmf}
%		\begin{align}
%			0&=\frac{\Delta}{\taum \pi}+2rv\\
%			0&=		v^2  + J\taum uxr - (\pi \taum r)^2 +I_1\\
%			0&=	(1-x)/\taud-uxr\\
%			0&=	(U_0-u)/\tauf+U_0(1-u)r
%		\end{align}
%	\end{subequations}
%	\htnote[inline]{write F(X) = 4 lines 12a-12d, same for G(X)}	
We can therefore already infer the shape of $S_0$.
It corresponds to the cartersian product $S^\ast \times \{I_2|I_2\in \mathbb{R}\}$, where $S^\ast$ denotes the S-shaped branch of equilibria of the fast subsystem and is given in \cref{eq:Sast}.
\begin{align}\label{eq:Sast}
	S^\ast:=\{(r,v,x,u,I_1)|\fastrhs(\Xf,\Xs)=0\}
\end{align}
%	\htnote[inline]{cross check $S^\ast$}
%	\htnote[inline]{change notation of $\dot{r}$ etc to F(X)=0}

\begin{figure}[h]
	\centering
	\includegraphics[width=1\linewidth]{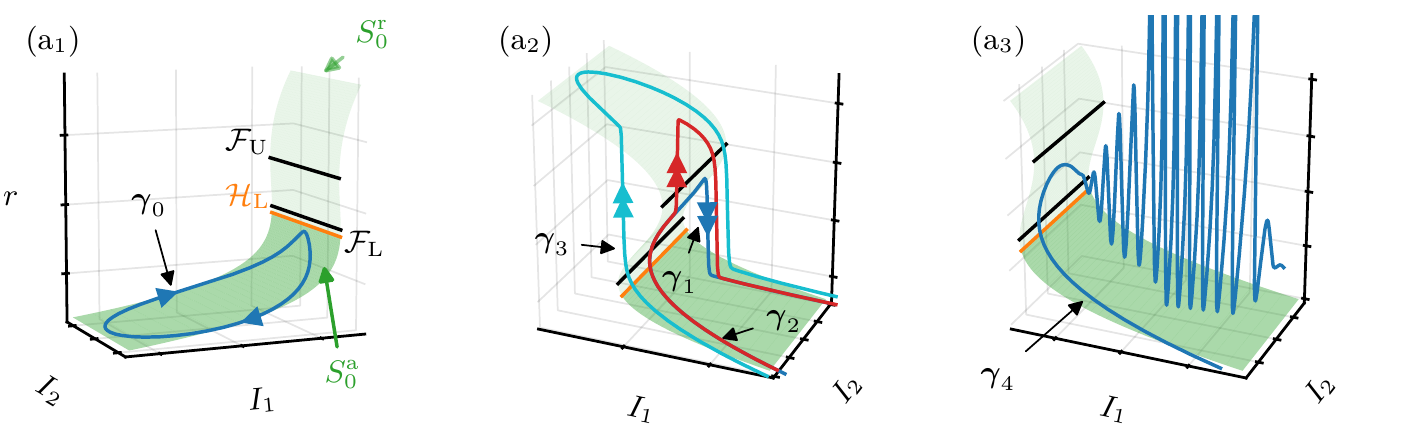}
	\caption{Slow-fast dissection and critical manifold: Solutions $\boldsymbol{\gamma}_0(t)$ to $\boldsymbol{\gamma}_4(t)$ of the full system superimposed on the critical manifold $\cmf$ in $(I_1,I_2,r)$-space.
	The parameter values are identical to those in \cref{fig:fullsystemsolutions1}.}
	\label{fig:fullsystemsolutions2}
\end{figure}

Hence the associated fold set $\mathcal{F}$ has two 1D connected components, namely the two lines $\FUc:=\{\FU\} \times \{I_2|I_2\in \mathbb{R}\}$ and $\FLc:=\{\FL\} \times \{I_2|I_2\in \mathbb{R}\}$. This means that the two fold lines are parametrised by $I_2$. These two fold lines of $\cmf$, along $\FLc$ and $\FUc$, can be seen in \cref{fig:fullsystemsolutions2}.

Points of $\cmf$ are equilibria of the fast subsystem. This mean that, in addition, local stability properties of the fast subsystem can be associated with points of the critical manifold. Stable (unstable) parts of the equilibrium branch in \cref{fig:bifdiag} will become attracting (repelling) sheets of the critical manifold. This property can be seen as an indicator for fast flow in the full system, distant from $\cmf$. Here fast dynamics act and $(r,v,x,u)$ will be repelled and attracted accordingly. The stability changes along the set of Hopf bifurcation given by $\HBLc:=\{\HBL\} \times \{I_2|I_2\in \mathbb{R}\}$ and $\HBUc:=\{\HBU\} \times \{I_2|I_2\in \mathbb{R}\}$ . Attracting (repelling) parts of $\cmf$ are marked as green (light green) surfaces throughout the present work. \\

\subsection{Singular dynamics: fast subsystem}
The extended model \cref{eq:model1} is able to generate periodic oscillations due to plastic synapses even in absence of time dependent forcing ($\Iext=\text{const.}$), as outlined in \cref{sec:neurmassSTP}. Their existence depends on the exact choice of parameters values, one of the important ones being the total nonsynaptic current given by $\bar{\eta}+\Iext$. Limit cycles can arise via a plethora of bifurcation scenarios. Here we considered the case of a subcritical Hopf bifurcation followed by a fold of limit cycles, giving rise to stable oscillatory behaviour, when considering $I_1$ as a bifurcation parameter (see \cref{fig:bifdiag}). Since $\cmf$ consists of fast subsystem equilibria, the bifurcation diagram is essentially a projection of $\cmf$ onto the $I_1$-$r$ plane.

\subsection{Singular dynamics: slow subsystem}

\subsubsection*{Slow flow.}
We make further use of the dissection by studying the \textit{slow flow} on the critical manifold $\cmf$. In the slow subsystem \cref{eq:slowsub} the state space is reduced to $\cmf$, described by four algebraic conditions in \cref{eq:critmf}, the solutions of which depend on the slow variable $I_1$ entering into the membrane potential equation. The state variables in this limit are subject to the slow flow $(\dotsl{\Xf},\dotsl{\Xs})$ describing their dynamics on $\cmf$. For $\Xs=(I_1,I_2)$ this is explicitly given via the Hopf normal form \cref{eq:forcing2}. For the fast variables $\Xf=(r,v,x,u)$ however, the algebraic constraints define $\Xf$ as well as $\dotsl{\Xf}$ on $\cmf$ implicitly. In this case the flow can be obtained by taking the total (slow) time derivative of \cref{eq:slowsubcritmf} as done in \cref{eq:slowflowderiv1},
\begin{align}
	0 =\dv{}{\tslow} \fastrhs(\Xf(\tslow),\Xs(\tslow)) =	\pdv{\fastrhs}{\Xf} \dv{\Xf}{\tslow} + \pdv{\fastrhs}{\Xs} \dv{\Xs}{\tslow} \quad,  \label{eq:slowflowderiv1}
\end{align}
where $\pdv{(\cdot)}{\boldsymbol{a}}$ is the Jacobian of $(\cdot)$ with respect to $\mathbf{a}$. If the Jacobian $\pdv{\fastrhs}{\Xf}$ is invertible, i.e., $\det (\pdv{\fastrhs}{\Xf})\neq 0$, then the slow flow of $\Xf$ can be calculated and results in \cref{eq:slowflowfastvar}.
\begin{align}
	\dv{\Xf}{\tslow}\equiv\dsXf = -\left(\pdv{\fastrhs}{\Xf}\right)^{-1} \left(\pdv{\fastrhs}{\Xs}\dsXs\right)\label{eq:slowflowfastvar}
\end{align}
This slow flow is only defined on $\cmf$ and represents a system of ODEs capturing the dynamics of $\Xf$ and $\Xs$ on the manifold.

In the particular case of the neural mass we have $\Xf=(r,v,x,u)$, $\Xs=(I_1,I_2)$ and $\fastrhs$, $\slowrhs$ given in \cref{eq:model2,eq:forcing2}, respectively.
Values of $\Xf$ on $\cmf$ are determined via the algebraic conditions \cref{eq:critmf}, which entangle the fast variables to each other and to $I_1$. Given either $r$, $v$, $x$ or $u$, the other components can be calculated straightforwardly. In other words, it is sufficient to consider the slow flow of one the fast variables, here $r$, to understand the slow dynamics. Taken this into account we finally obtain the slow flow $(\dotsl{r},\dotsl{I_1},\dotsl{I_2})$ given in \cref{eq:sflow}
\begin{subequations}\label[pluralequation]{eq:sflow}
	\begin{align}
		\dotsl{r}   & =g_1(I_1,I_2)Ar(1+ru\taud)(1+rU_0\tauf)/D \\
		%			\dotsl{v}&=-g_1(I_1,I_2)A I_2 v(1+ru\taud)(1+rU\tauf)/D\\
		%			\dotsl{d}&=g_1(I_1,I_2)A I_2 dr\taud(u+rU\tauf)/D\\
		%			\dotsl{u}&=g_1(I_1,I_2)A I_2 Ur\tauf(1-u)(u+ru\taud)/D\\
		\dotsl{I_1} & =g_1(I_1,I_2)                             \\
		\dotsl{I_2} & =g_2(I_1,I_2)
	\end{align}
\end{subequations}

We can see the relevance of the denominator $D$, given in \cref{eq:sflowdenom}, by noting that $D=0$ defines the \textit{singular points} of $\cmf$ at which the fast subsystem Jacobian $\pdv{\fastrhs}{\Xf}$ is singular.
\begin{align}
	D=2 \left[(\pi  r)^2+v^2\right] (r \taud u+1) (r \tauf U_0+1)-Jx  r (r \tauf U_0+u)\label{eq:sflowdenom}
\end{align}
Saddle-node bifurcations of the fast subsystem are characterized by the identical condition $\det (\pdv{\fastrhs}{\Xf})=0$ and singular points are equivalent to the fold curves $\FLc$ and $\FUc$ shown in \cref{fig:fullsystemsolutions2}. Therefore, the slow flow is undefined along these lines and the slow subsystem fails to describe the slow dynamics for trajectories intersecting with $\FLc$ or $\FUc$.

\subsubsection*{Desingularization.}\label{sec:desingularization}
This limitation can be mitigated by introducing an auxiliary system and desingularizing \cref{eq:sflow} through the application of a nonlinear time rescaling $\tau \mapsto D\cdot\tau$. One obtains the \textit{desingularized reduced system} (DRS) given in \cref{eq:drsflow} with $\hat{\tau}:=D\tau$.
\begin{subequations}\label[pluralequation]{eq:drsflow}
	\begin{align}
		\frac{\dd r}{d\hat{\tau}}   & =g_1(I_1,I_2)A r(1+ru\taud)(1+rU_0\tauf) \\
		%			\frac{\dd v}{d\hat{\tau}}&=-g_1(I_1,I_2)A v(1+ru\taud)(1+rU\tauf)\\
		%			\frac{\dd x}{d\hat{\tau}}&=g_1(I_1,I_2)A dr\taud(u+rU\tauf)\\
		%			\frac{\dd u}{d\hat{\tau}}&=g_1(I_1,I_2)A Ur\tauf(1-u)(u+ru\taud)\\
		\frac{\dd I_1}{d\hat{\tau}} & =g_1(I_1,I_2)\cdot D                     \\
		\frac{\dd I_2}{d\hat{\tau}} & =g_2(I_1,I_2)\cdot D
	\end{align}
\end{subequations}
The DRS benefits from the fact that the singularities are resolved, allowing to investigate the slow dynamics near and on the fold lines $\FLc$ and $\FUc$. At the same time new equilibria are introduced satisfying $D=0$.
Additionally, as a consequence of the employed non-linear time rescaling, the flow direction is not preserved. At the fold curves, with $D=0$, a change of sign of $D$ takes place. Hence, between $\FLc$ and $\FUc$, i.e, on the middle sheet of $\cmf$, the flow of the DRS is opposite to the slow flow.

Slow trajectories which entirely remain on the same sheet of $\cmf$ can be easily understood using the slow flow \cref{eq:sflow}. Canard orbits evolve on attracting as well as repelling sheets of $\cmf$ and therefore require a view on the dynamics near the folds. For this, we will determine the equilibria of the DRS \cref{eq:drsflow} in the following and analyze their invariant manifolds. Equilibria of the DRS which satisfy $D=0$ necessarily coincide with a fold of $\cmf$. As we will show in the following section, this gives rise to so-called \textit{folded singularities}.

\subsubsection*{Folded saddle and folded homoclinics.}\label{sec:foldedsaddlehom}
Up to three focus equilibria are located at $p_1=(I_1,I_2,r)=(0,0,r_k)$, where the $r_k$ are points of $\cmf$ given $(I_1,I_2)=(0,0)$. In this work we will remain at  $\bar{\eta}$ values for which only one equilibrium $p_0=(0,0,r_0)$ exists. This point is the only equilibrium of the slow flow \cref{eq:sflow}. In the DRS (\cref{eq:drsflow}) however, the condition $\{g_1(I_1,I_2)=0,D=0\}$ yields additional fixed points $p_1$ and $p_2$ located on the fold lines at $D=0$.

\begin{figure}[h]
	\centering
	\includegraphics[width=1\linewidth]{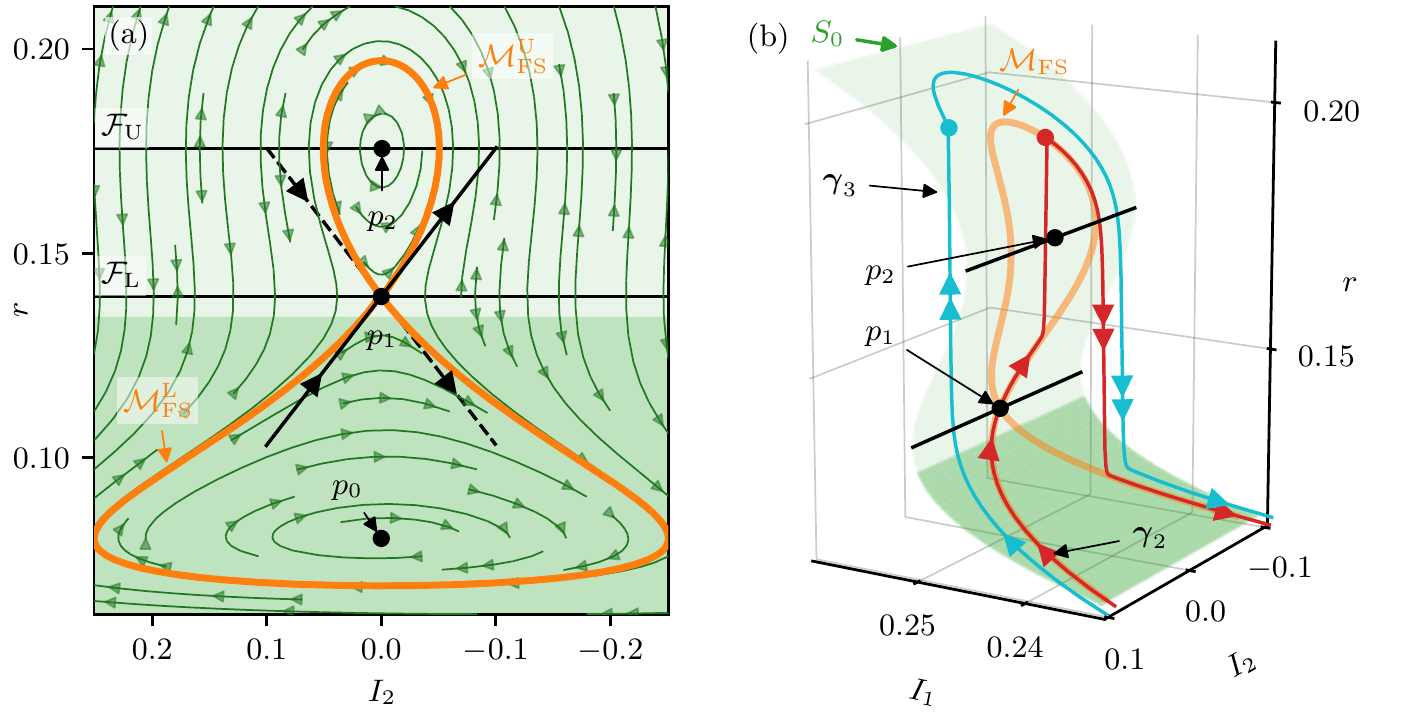}
	\caption{
	Folded-saddle and jump-on canards:
	{\bf (a)} Critical manifold $\cmf$ in $r-I_2$ projection superimposed with the slow flow [green arrows, see \cref{eq:sflow}]. The black dots $p_0$ (unstable focus), $p_1$ (saddle), $p_2$ (center) denote equilibria of the DRS [\cref{eq:drsflow}]. The point $p_1$ denotes a folded-saddle equilibrium with the associated stable (unstable) eigendirection indicated by a solid (dashed) arrow along the slow flow. The orange curves $\FSmf$ mark the stable and unstable manifolds of $p_1$, forming heteroclinic connections through $p_1$.
	{\bf (b)} $\cmf$ in $(I_1,I_2,r)$-space. The curves $\boldsymbol{\gamma}_2(t)$ and $\boldsymbol{\gamma}_3(t)$ are solutions of the full system.
	The objects $p_1,p_2$ and $\FSmf$ depend on the choice of $A$, here they correspond to the value used to obtain $\boldsymbol{\gamma}_3(t)$. Other parameters values are as in \cref{fig:fullsystemsolutions1,fig:fullsystemsolutions2}.
	}
	\label{fig:DRSFoldedSaddle}
\end{figure}

The three equilibria of the DRS are displayed in \cref{fig:DRSFoldedSaddle}(a) on $\cmf$ in $r-I_2$ projection. In the full system, for sufficiently small $A$ and $\eps$, solutions lie close to the bottom, attracting, sheet of $\cmf$. When the amplitude is increased, these cycles can pass very close to $p_1$ and start to follow the middle, repelling sheet. One way to understand this canard dynamics is to make use of the properties of the $p_k$ in the DRS and their role for the slow subsystem.

Located on the bottom sheet of $\cmf$, $p_0$ results from the Hopf form given in \cref{eq:forcing1} and is an unstable focus at ($I_1,I_2)=(0,0)$. On the other hand, $p_2$ lies on the upper fold line $\FUc$ and denotes a center, i.e., it has purely imaginary complex conjugate eigenvalues. The equilibrium $p_1$ can be found on the lower fold line $\FLc$ and is of saddle type.
At a specific value of $A$, namely, when the forcing cycle intersects with $p_1$, an 8-shaped double homoclinic connection $\FSmf=\FSmf^{\rm L}\cup \FSmf^{\rm U}$ forms, consisting of two parts, which are connected via $p_1$; see the orange curve in \cref{fig:DRSFoldedSaddle}(a,b). The connection $\FSmf^{\rm L}$ is located on the lower sheet, while  $\FSmf^{\rm U}$ spans the middle and upper sheet of $\cmf$. They revolve around the unstable focus $p_0$ and center $p_2$, respectively, and are the stable and unstable manifolds of $p_1$.

The points $p_2$, $p_1$, and in particular the invariant manifolds associated with $p_1$, play an important role for the slow subsystem. Due to the negative sign $D<0$ on the middle sheet of $\cmf$ the slow flow is reversed with respect to the DRS. As a consequence the DRS saddle $p_1$ and center $p_2$ become folded singularities of the slow subsystem. These \textit{folded saddle} ($p_1$) and \textit{folded center} ($p_2$) are not equilibria of the slow flow.	However, for the slow dynamics they have similar impact on the dynamics as there unfolded counterparts, but with the crucial difference of reversed flow direction between $\FLc$ and $\FUc$. Accordingly, the folded saddle $p_1$ has significant influence on the dynamics of the slow subsystem along $\FSmf$, as described in the following.

(i) First of all, trajectories in the DRS evolving on $\FSmf$ necessarily approach $p_1$ asymptotically from the direction of the stable eigenvector, but can never pass through the saddle.

(ii) In the slow subsystem however, the folded saddle $p_1$ allows a pass-through along this direction. Below $\FLc$ trajectories are attracted to and above repelled from $p_1$.

(iii) The double homoclinic connection $\FSmf$ of the DRS is referred to as a \textit{folded homoclinic} in the slow subsystem. For this solution the passage of trajectories through $p_1$ occurs in finite time \cite{desrochesSpikeaddingParabolicBursters2016}.

(iv) Using the same type of argument, the invariant manifold $\FSmf^{\rm U}$ around the folded center $p_2$ becomes disconnected at the two intersections with $\FUc$, due to a reversal of the slow flow direction.  Solutions of the slow subsystem on $\FSmf^{\rm U}$ can not cross this line.

As a consequence of the previous properties (i)-(iv),  a singular canard exists in the slow subsystem, given a specific value of $A$. It evolves along the folded homoclinic $\FSmf$ below the bottom fold line $\FLc$ and extends, while remaining on $\cmf$, beyond the folded-saddle $p_1$ until the upper fold $\FUc$.

\subsubsection*{Singular canard orbits.}\label{sec:singularcanardorbits}
For the construction of singular orbits, we note once more that the middle sheet of $\cmf$ is repelling while the bottom sheet is attracting. Accordingly, a continuum of fast segments emerging from the middle sheet and connecting to the bottom sheet exist in the fast subsystem. This family of fast orbits collides with $\FSmf$ and likewise with the singular canard described above. As a result, infinitely many singular orbits can be constructed, by merging the singular canard at arbitrary positions on the middle sheet of $\cmf$, with fast segments.

These singular orbits evolve on the bottom sheet of $\cmf$, continue through the folded-saddle $p_1$,  while following the folded homoclinic, and jump at different heights, in terms of the coordinate $r$, from the middle to the bottom sheet of $\cmf$.
The full system solution $\boldsymbol{\gamma}_1(t)$ in \cref{fig:fullsystemsolutions1,fig:fullsystemsolutions2} displays this type of dynamics. The singular canard, hence also the family of singular canard orbits, can at most reach the upper fold line $\FU$. Here the slow flow is undefined and the reduction of state space to $\cmf$ fails to describe the dynamics. This is additionally reflected by the fact that $\FSmf^{\rm U}$ is disconnected at the intersections with $\FUc$, due to the folded property of $p_2$; see (iv) above. The singular canard orbit which reaches up until this point is the maximal canard.

\section{Full system dynamics: beyond singular orbits and classical canards}\label{sec:fsdynamics}

%\subsection{Beyond singular orbits}\label{sec:fenichel}
The solutions $\boldsymbol{\gamma}_0(t)$ to $\boldsymbol{\gamma}_4(t)$ are results of numerical computations for $\eps>0$.  As such, their slow segments evolve not on, but in an $\mathcal{O}(\eps)$ neighborhood of $\cmf$. This is to some extent an implication of Fenichel's theory \cite{fenichelGeometricSingularPerturbation1979}. For $0<\eps\ll 1$, it guarantees the existence of a \textit{slow manifold} $S_\eps$, that is in an $\mathcal{O}(\eps)$ neighborhood of $\cmf$, if $\cmf$ is \textit{normally hyperbolic} (see below).  Additionally, $S_\eps$ is locally invariant under the flow of the full system. The theorem also states that stable and unstable manifolds associated to $\cmf$ persist as $\mathcal{O}(\eps)$ perturbations. In other words, the flow on $S_\eps$ can be seen as a perturbation of the flow on $\cmf$; and the flow perpendicular to $\cmf$ as a perturbation of the fast subsystem's flow. For normally hyperbolic critical manifolds, one can deduce that singular orbits persist for $\eps$ and perturb into an $\mathcal{O}(\eps)$ neighborhood.

Normal hyperbolicity requires all eigenvalues of the Jacobian $\pdv{\fastrhs}{\Xf}|_{\cmf}$ to have non-zero real part \cite{hekGeometricSingularPerturbation2010}. The theorem can therefore not be applied on $\FLc$ and $\FUc$, given that they describe lines of saddle-node bifurcations. However, one can consider the three sheets of $\cmf$ separately, each one up to an $\mathcal{O}(\eps)$ neighborhood of the folds. From this we can conclude the persistence of slow segments nearby $\cmf$ for $\eps>0$, including the repelling segments within a canard solution, until close to the folds. Fenichel's theory does not encompass whether or not a connection of these segments exists. Yet, $\boldsymbol{\gamma}_1(t)$ to $\boldsymbol{\gamma}_3(t)$ exemplify such connection numerically: segments that evolve close to the bottom sheet of $\cmf$ connect to segments close the middle sheet. A detailed treatment of these connections nearby the non-hyperbolic points, exceeds the scope of this work. For more rigorous approaches, we refer to non-standard analysis \cite{benoitChasseAuCanard1981}, matched asymptotics \cite{eckhausRelaxationOscillationsIncluding1983} and so-called blow-up techniques \cite{krupaExtendingGeometricSingular2001}. With these advanced methods, it is possible to show that orbits with canard segments of different length perturb at different parameter values within an exponentially narrow regime, thus leading to the canard explosion in the full system.

\subsection{Jump-on canards}
The construction of singular canard orbits linked with Fenichel's theorem explains the dynamics of $\boldsymbol{\gamma}_1(t)$ in \cref{fig:fullsystemsolutions1,fig:fullsystemsolutions2}, where the forcing amplitude $A$ is large enough to surpass the lower Hopf bifurcation $\HBL$ and the lower fold $\FL$. \textit{Headless} canards, like this one, have a jump to the bottom branch in common.
Slow-fast systems may also have canard solutions with a head.  They usually appear in systems, which have two folds: the first one destabilizing, the second one stabilizing the branch. Headed canards jump onto this stable upper part. In 3D systems with an 1D or 2D S-shaped critical manifold, the upper sheet is typically unstable near the upper fold, thus preventing the existence of headed canards. Instead, past the maximal canard, fast oscillations related to the existence of limit cycles develop and lead to bursting solutions; see \cref{sec:transbursting}.

The trajectory $\boldsymbol{\gamma}_2(t)$ in \cref{fig:fullsystemsolutions1,fig:fullsystemsolutions2} has the characteristic dynamics of a headed canard. However, it is peculiar for various reasons. First of all, we can classify this type of dynamics as \textit{jump-on canard}, since the trajectory lands (after the fast jump) on a seemingly repelling slow manifold. Moreoever, the jump-on dynamics can occur after a regular canard segment, as for $\boldsymbol{\gamma}_2$ or independent of that, like for $\boldsymbol{\gamma}_3$. As a matter of fact, trajectories of the latter type resemble relaxation oscillations, like in VdP, despite the repelling upper sheet of $\cmf$. This has an additional consequence: for small enough $\eps$ a continuous transition from subthreshold oscillations to bursting is blocked by jump-on canards.
The fact that fast oscillations other than relaxation oscillations are absent beyond the maximal canard are novel and unexpected phenomena.
In the following we will address how these solutions emerge. Their impact on the route towards bursting is discussed in details in \cref{sec:blocking}.

%	To our knowledge canards connecting two repelling sheets of the critical manifold, as well as the absence of fast oscillations
%	\mdnote{we need to be a bit more precise here: the absence of fast oscillations other than a relaxation oscillation, which is what we get with these jump-on canard solutions} beyond the maximal canard in this system are novel and unexpected phenomena.

The two solutions $\boldsymbol{\gamma}_2(t)$ and $\boldsymbol{\gamma}_3(t)$ are shown in \cref{fig:DRSFoldedSaddle}(b) and \cref{fig:DRSJumpOn} using two different projections. They exemplify two types of jump-on canards, that have an approach towards globally repelling equilibria of the fast subsystem in common. The cyan solution $\boldsymbol{\gamma}_3$ represents an orbit which does not interact with the folded-saddle $p_1$. It slowly evolves on the lower sheet of $\cmf$, crosses the curve of Hopf bifurcations $\HBLc$ and reaches the lower fold curve $\FLc$ at which the slow subsystem is singular. Fast dynamics come into play and expectedly the dynamics will approach attractors of the fast subsystem. Such attractors for the considered $I_1$ value at which the curve escapes $\cmf$ are solely the stable limit cycles displayed in \cref{fig:bifdiag}(a).
Instead of entering a period of bursting, the cyan trajectories approaches the upper branch of unstable equilibria. As soon as it jumps onto on $\cmf$, the slow subsystem becomes a valid limit anew and the curve remains on $\cmf$ until it reaches $\FUc$, where it jumps down to the stable sheet.

In the second case, $\boldsymbol{\gamma}_2$ in \cref{fig:DRSFoldedSaddle}(b,c), the orbit possesses a canard segment and jumps from the repelling middle sheet to the repelling upper sheet of $\cmf$. Similar to the previous case, it evolves on $\cmf$ until $\FUc$ and finally jumps down. The global motion is identical to that of a canard with a head in VdP.

Multiple elements of these singular cycles have to be understood.
First of all, both cases have a similar slow segment in common, namely the part of the trajectory on the upper sheet of $\cmf$. They can be approximated by solutions of the reduced problem \cref{eq:slowsub} and are enforced by the presence of $p_2$, around which the trajectory evolves. Since the center $p_2$ is folded, full rotations around it are not possible and the slow parts terminate at $\FUc$, where the slow flow is undefined. Here the trajectory can be joined to a fast bit which connects from $\FUc$ to the attracting sheet of $\cmf$. After this part, the dynamics on $\cmf$ is again governed by the slow subsystem and depending on the forcing amplitude $A$, the two orbits take different paths. The solution $\boldsymbol{\gamma}_3$ crosses $\FLc$ far from $p_1$ and the slow segment on the attracting sheet stops; $\boldsymbol{\gamma}_2$ passes through a neighbourhood of $p_1$ and exhibits canard dynamics before a jump occurs. These parts of the orbits are entirely described within the scope of the slow subsystem.

\subsection{Nested timescale separation}\label{sec:jump-on}
\begin{figure}[h]
	\centering
	\includegraphics[width=0.5\linewidth]{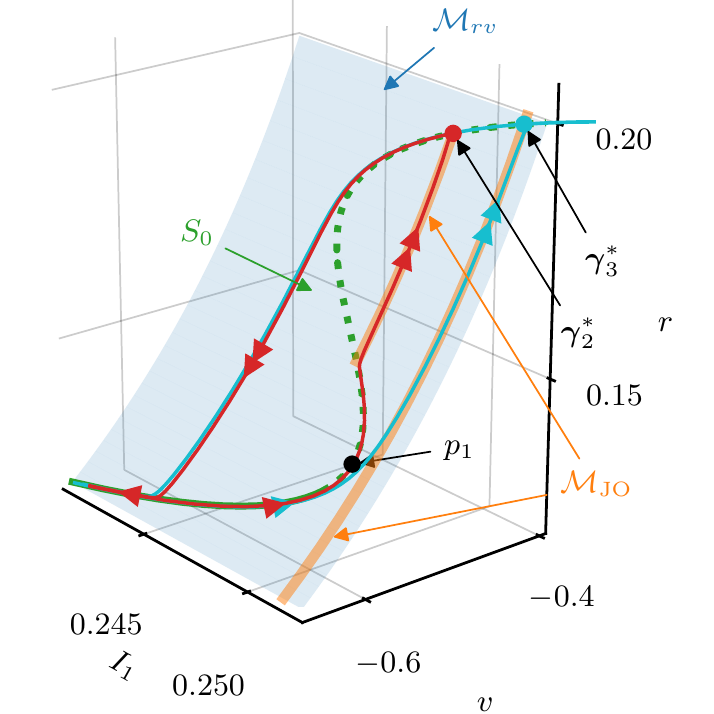}
	\caption{
		Folded-saddle and jump-on canards:
		$\cmf$ in $I_1$-$v$-$r$ space (green curve). The blue surface $\rvmf$ is defined by $r=-\frac{\Delta}{\pi \taum v}$. The curves $\boldsymbol{\gamma}_{2,3}(t)$ are solutions of the full system and their jump-on points $\boldsymbol{\gamma}_{2,3}^\ast$ on the upper repelling sheet of $\cmf$ are marked by dots. Parameters values are as in \cref{fig:DRSFoldedSaddle}. 
	}
	\label{fig:DRSJumpOn}
\end{figure}

Understanding the remaining segment that leads $\boldsymbol{\gamma}_2$ and $\boldsymbol{\gamma}_3$ towards the repelling sheet of $\cmf$ requires a more detailed analysis. The mechanism is the same for both cases and will be discussed in the following. Since these pieces of the orbits evolve on the fast timescale we present a visualization of the curves in $I_1-v-r$ space in \cref{fig:DRSJumpOn}. In this projection the critical manifold $\cmf$ is shown as a green curve $r(I_1), v(I_1)$ with attracting (repelling) parts as a solid (dashed) line. In the singular limit, the jump-on points $\boldsymbol{\gamma}^\ast_{2,3}=(r_{2,3}^\ast,v_{2,3}^\ast, x_{2,3}^\ast, u_{2,3}^\ast)$  of the jump-on canard solutions $\boldsymbol{\gamma}_2$ and $\boldsymbol{\gamma}_3$, marked by the red and cyan dots in \cref{fig:DRSJumpOn}, are of saddle-focus type. Linearization of the dynamics reveals a weakly and strongly attracting direction in a neighbourhood of the jump-on points together with a repelling direction with complex conjugate eigenvalues. This suggests a 2D stable manifold leading to $\boldsymbol{\gamma}^\ast_{2,3}$.
We will simplify the problem further by noting that the entire dynamics of jump-on canards takes place close to the surface $\rvmf$ defined in \cref{eq:defMrv}. On the one hand, for the slow pieces of the curve this observation is as expected, since by definition this conditions holds on $\cmf$.

\begin{align}
	\dotfa{r}=0 \Rightarrow r=-\frac{\Delta}{2\pi \taum v} \text{\quad for } v\neq 0\label{eq:defMrv}
\end{align}
On the other hand, also the fast parts of the orbits remain on $\rvmf$. This implies a reduction of the fast dynamics to $\rvmf$ for the part of state space in which jump-on canards can be found. We will exploit this reduction and investigate a secondary differential-algebraic system resulting from the fast subsystem \cref{eq:fastsub} in which the flow $\dotfa{r}$ is equilibrated.

This latter step, without having a complete picture of the time scaling, is analogous to an additional dissection of the fast subsystem. In other words, the full system exhibits three timescales for what concerns jump-on canards: the dynamics of $r$ takes place on a fast, that of $(v,x,u)$ on an intermediate and that of $(I_1,I_2)$ on a slow timescale. The equilibration $\dotfa{r}=0$ eliminates the fastest of these scales and approximates the intermediate scale dynamics of jump-on canards, which takes place in the vicinity of $\rvmf$.

In this new framework $\rvmf$ describes a manifold on which the dynamics of $(v,x,u)$ take place, while $I_1, I_2$ remain frozen.
The jump-on points $(v^\ast, x^\ast, u^\ast)$ [red and cyan dots in \cref{fig:DRSJumpOn}] are mutual points of $\rvmf$ and the upper sheet of $\cmf$. In the reduced problem on $\rvmf$,  $(v^\ast, x^\ast, u^\ast)$ are of saddle type with eigenvalues $\lambda_1 \gg \lambda_2\gg -\lambda_3 > 0$ and therefore have associated 1D stable manifolds $\JOmf\subset \rvmf$ [orange curves in \cref{fig:DRSJumpOn}]. They exists for any value of $I_1$ after the upper fold. For values of $I_1$ beyond the lower fold, they extend down to $r=0$. The family of 1D stable manifolds associated with the jump-on points can guide trajectories towards the upper repelling sheet of $\cmf$ and this way leads to the existence of jump-on canards.

So far we have discussed singular orbits, for which the fast segments connect different sheets of the critical manifold $\cmf$. In the case of regular canards, for $\eps>0$ small enough, these correspond to stable equilibria of the fast subsystem. For jump-on canards they might be unstable, but possess a stable direction, allowing to reach and stay on the repelling sheet of $\cmf$. The main dynamics of these singular orbits takes place on $\cmf$ and does not display phases of fast oscillations, as for bursting solutions, but only single fast jumps.

\section{Slow-fast transition to bursting: a tale of two routes}\label{sec:transbursting}
Opposed to that, the solution $\boldsymbol{\gamma}_4(t)$ [see \cref{fig:fullsystemsolutions1,fig:fullsystemsolutions2}] and the case initially illustrated in \cref{fig:burst_intro}  exhibit bursts: a slow segment is followed by fast oscillations. Periodic solutions of the fast subsystem are the underpinning elements of bursting, such that a classification in terms of the fast subsystem's bifurcations appears appropriate. As shown in the bifurcation diagram \cref{fig:bifdiag}(a), limit cycles originate and terminate at Hopf bifurcations, and change stability at a fold of cycles. Strictly following the classification of Izhikevich, bursting solutions in this system are of subcritical Hopf/fold cycle type. However, the subcritical Hopf bifurcation is closely followed by a fold of the underlying equilibrium branch. Due to a delay effect when surpassing the subcritical Hopf bifurcation \cite{neishtadtProlongationLossStability1987,neishtadtProlongationLossStability1988,baerSlowPassageHopf1989}, bursting can effectively be initiated at that fold; see e.g. \cref{fig:spike-adding_1em3}(b$_4$). We will restrict our analysis to these cases.
Here a more aptly description of the bursting type is fold/fold cycle, which corresponds to elliptic bursting in the classification of Rinzel \cite{rinzelFormalClassificationBursting1987}.

In order to understand bursting solutions of the full system, we want to remain in the slow-fast dissection. However, the neural mass with STP in presence of periodic forcing turns out to be a peculiar system and numerically difficult to handle. We are constrained by two main factors. First of all, as soon as we leave the singular limit, i.e, for $\eps > 0$, slow segments of trajectories diverge sensitively from the critical manifold. In other words, $\eps$ is required to be remarkably small to maintain a good agreement between between full system trajectories and singular orbits. Secondly, numerical simulation as well as numerical continuation of the full system for small enough $\eps$ are challenging, since the dynamics appears to be stiff and require high accuracy.

Therefore a clear view on the emergence of bursting can not be gained easily in this framework. For this reason we will provide, additionally to geometrical arguments, numerical evidence on how bursting forms in the present system, either via direct simulation or continuation using the full system.
In general, bursts might emerge via a \textit{spike-adding} mechanism, that is, the consecutive addition of spikes into the orbit, when varying a parameter (e.g., the forcing amplitude).	For parabolic bursting this spike-adding is mediated by folded-saddle canards \cite{desrochesSpikeaddingParabolicBursters2016}; in square-wave bursters on the other hand, passages through a fast-subsystem saddle-homoclinic bifurcation and a folded node determine the number of large-amplitude oscillations in the burst and small-amplitude oscillations before the burst, respectively \cite{desrochesMixedmodeBurstingOscillations2013}.
In the following, we report the spike-adding mechanism for the NMSTP. At its basis is an interaction of the canards dynamics invoked by the presence of the folded saddle $p_1$, as well as unexpected torus-canard dynamics. Moreover, we will point out the role of jump-on canards for this spike-adding transition.

\subsection{Canard explosion and spike-adding}\label{sec:bursting}

\begin{figure}[h!]
	%		\begin{adjustwidth}{-1.4in}{-0in}
	\centering
	\includegraphics[width=1\linewidth]{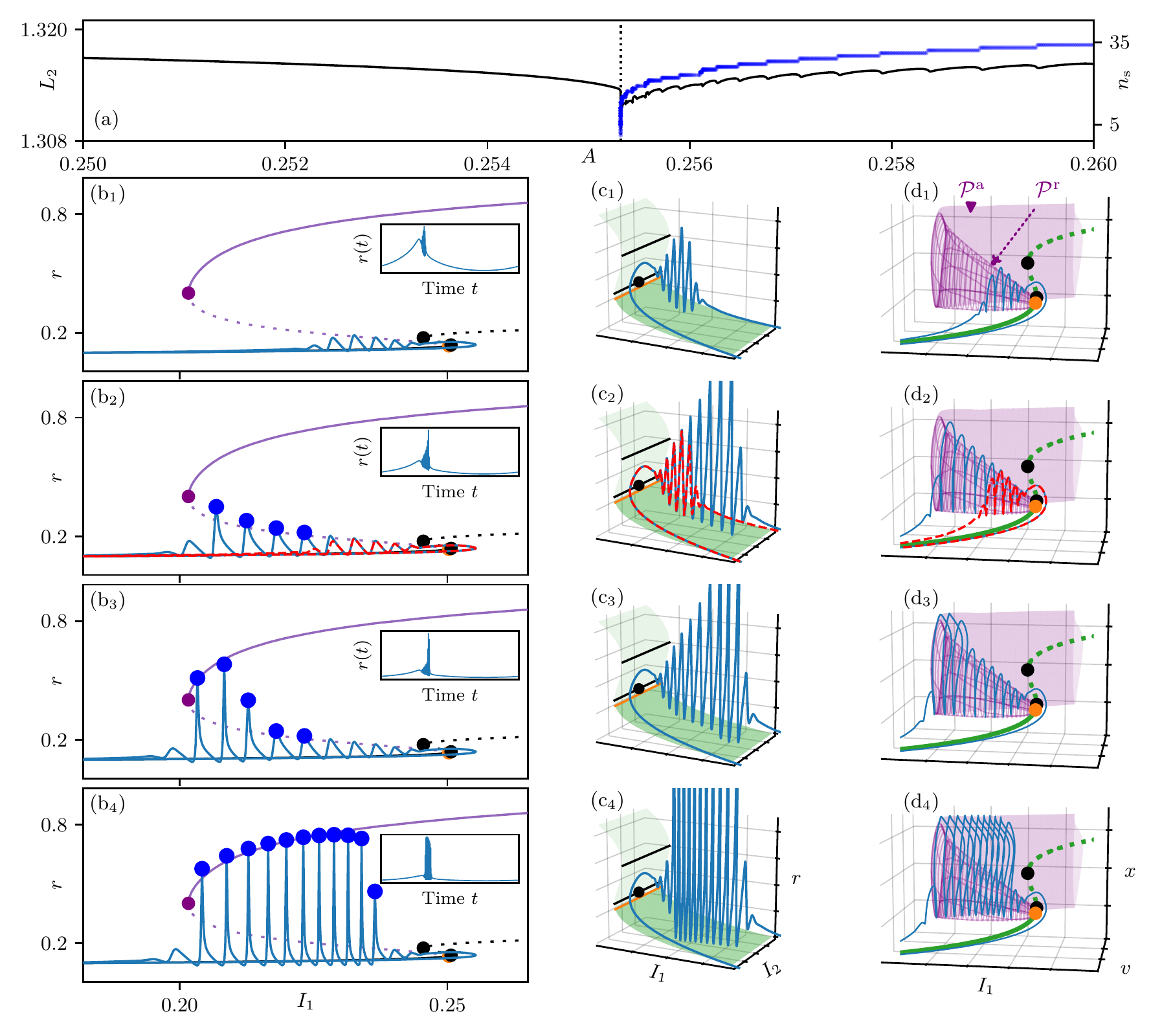}
	\caption{Emergence of bursting: {
	\bf (a)} Bifurcation diagram of the full system. The black curve shows the $L_2$-norm of a family of periodic solutions vs. the forcing amplitude $A$. The blue dots show the number of spikes $\nspikes$, defined as the local maxima of for which $r(t)>\SI{0.21}{\taum^{-1}}$. The dashed vertical line is located at $A=A^\ast=0.25531851205$. {\bf (b)} Solutions $(r(t),v(t),x(t),v(t),I_1(t),I_2(t))$ in $r-I$ projection superimposed on the bifurcation diagram of the fast subsystem. The insets show the solution in time.
		{\bf (c - d)} Same solutions as in (b) and critical manifold $\cmf$ in $(I_1,I_2,r)$-space (c) and in $(I_1,v,x)$-space (d).
	In (d) attracting (repelling) sheets of $\cmf$ are visualized as a green solid (dotted) line; the purple surface (wireframe) represents the family of stable (unstable) limit cycles of the fast subsystem [purple branch in column (b)]. Note that here the $x$-axis is inverted.
	In (b) the black dots denote $\FLc$ and $\FUc$, the orange dot $\HBLc$, while the black dots in (c) show the folded singularity $p_1$, assuming $A=A^\ast$. Spikes contributing to $\nspikes$ are marked with blue dots. The dashed red curve shows the solution of the panels above.
	The $A$ values in (b - d) are in increasing order from top to bottom and near $A^\ast$. Full system solutions obtained at $\eps=10^{-3}$.
	}
	\label{fig:spike-adding_1em3}
	%		\end{adjustwidth}
\end{figure}

To start with, we consider the case $\eps=10^{-3}$ and investigate the full system dynamics by performing continuation with the forcing amplitude $A$ as a parameter. We want to stress that this $\eps$ value, although very small, proves to be rather distant from the singular limit and slow-fast dissection arguments have to be taken with caution. The initial solution is for an $A$ value corresponding to subthreshold oscillations, like $\boldsymbol{\gamma}_0(t)$ in \cref{fig:fullsystemsolutions1,fig:fullsystemsolutions2}, and is continued towards larger amplitudes.

As a solution measure the $L_2$-norm of this family is plotted vs. $A$ in \cref{fig:spike-adding_1em3}(a). The first part until $A=A^\ast\approx 0.25531851205$ is in the subthreshold regime. Around $A^\ast$ a very sharp transition occurs, resembling a canard explosion. In this transition region the orbits already exhibit first spikes, here defined as the number $\nspikes$ of local maxima of $r(t)$ for which $r(t)>0.21$. This is followed up by a series of arches (on the solution branch) at $A>A^\ast$.

The arches are clearly related to the addition of new spikes to the orbit: with every termination of an arch,  by means of a vertical dip of the curve, the number of spikes $\nspikes$ increases. This behaviour can be better observed for larger $A$, for which each arch is related to the adding of exactly one spike. Prior, the arches lay more dense along $A$ and the spike adding appears to be of more complex nature. Despite the fact that the points of $\nspikes$ vs. $A$ depend on the choice of the $r$-threshold for which a spike is counted, it is evident that spikes are added consecutively. It is also clear that these bursting solutions emerge, in a continuous manner, from subthreshold oscillations.

A more detailed view on the full system dynamics near the explosive transition around $A=A^\ast$ is given in the columns (b - d) of \cref{fig:spike-adding_1em3}. Column (b) shows the solution in $r-I_1$ projection superimposed on the fast subsystem's bifurcation diagram; column (c) in $(I_1,I_2,r)$-space together with $\cmf$; column (d) in $(I_1,v,x)$-space.

In projection onto the ($I_1,v,x$)-space, the critical manifold appears as a curve ($I_1, v(I_1),x(I_1))$). Additionally we show the family of fast subsystem limit cycles, which emerge at the lower subcritical Hopf bifurcation $\HBL$. Unstable periodic solutions of this branch will be denoted by $\Gamma^{\rm r}$, stable ones after the fold of cycles by $\Gamma^{\rm a}$. Embedding these solutions $\Gamma^{\rm a, r}(I_1,t)=(r,v,x,u)(I_1, t)$ into the state space of the full system one obtains the surface $\LCmf=\LCmfa \cup \LCmfr$, which consists of attracting and repelling parts $\mathcal{P}^{\rm a,r} = \{(\Gamma^{\rm a,r}(I_1,t), I_1) | t\in [0,T(I_1)]\} \times \{I_2|I_2\in \mathbb{R}\}$, corresponding to stable and unstable branches of the solution family, respectively. The period $T(I_1)$ of these cycles depends on $I_1$.

\subsubsection*{Onset of fast oscillations $\boldsymbol{({\rm b}_1, {\rm c}_1,{\rm d}_1)}$.} At the smallest of the chosen $A$ values near $A^\ast$ one can already observe fast oscillations, consisting of five not fully developed spikes. They occur after the trajectory has turned around the folded saddle $p_1$, marked by a black dot in $({\rm c}_1)$. It is this motion around $p_1$, taking place in the vicinity of the critical manifold $\cmf$, which has signs of a turning point, guiding the trajectory along the repelling sheet of $\cmf$. Taking the rather large $\eps$ value into account, this turn hints at a canard segment arising due to the presence of the folded saddle $p_1$. After this segment, in $(I_1,I_2,r)$-space (panel ${\rm c}_1$), the trajectory pierces through the repelling sheet of $\cmf$ and fast oscillations set in. These results suggest that bursting is initiated at the termination of a canard segment, closely following the repelling middle sheet of $\cmf$.

\begin{figure}[h!]
	\centering
	\includegraphics[width=1\linewidth]{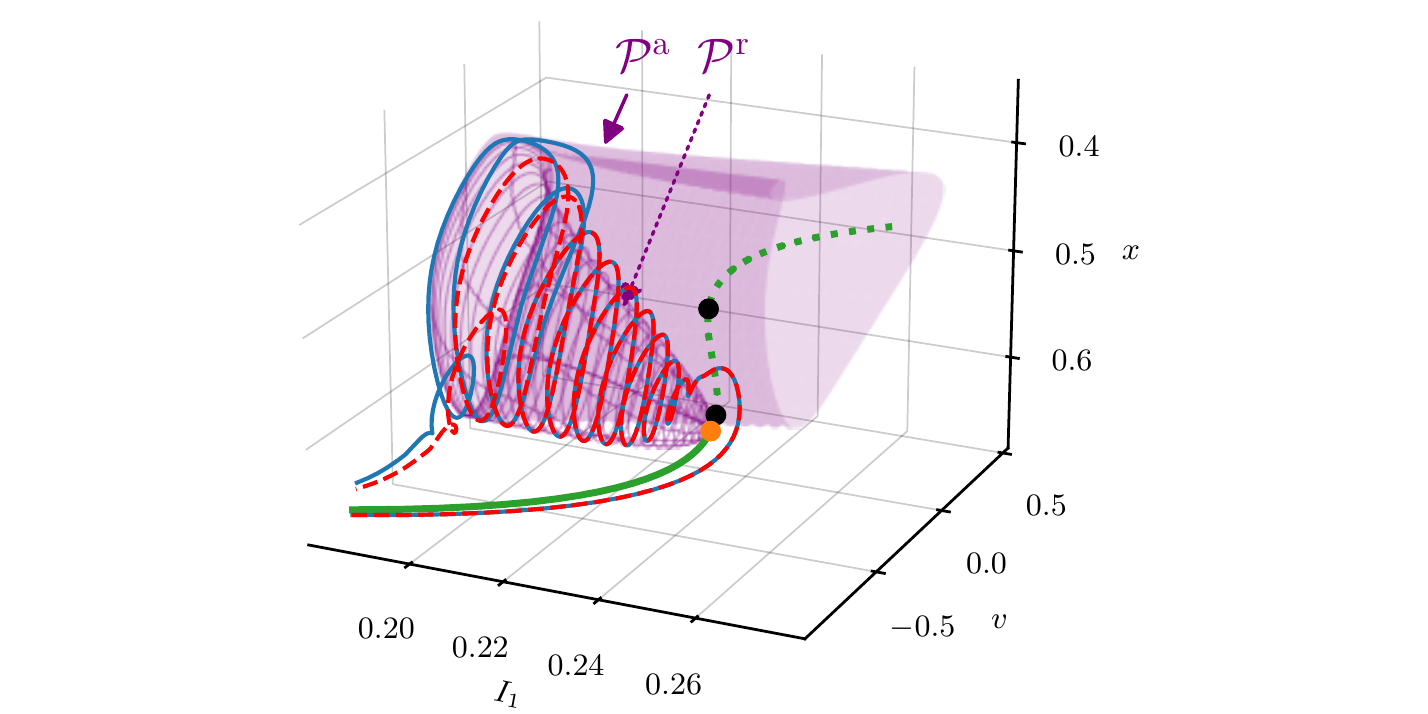}
	\caption{Emergence of bursting:
	Enlargement of \cref{fig:spike-adding_1em3}(${\rm d}_2$). Critical manifold in $(I_1,v,x)$-space with attracting (repelling) sheets of $\cmf$ visualized as a green solid (dotted) line; the purple surface (wireframe) represents the family of stable (unstable) limit cycles of the fast subsystem [purple branch in column (b)], emerging from the lower Hopf bifurcation $\HBL$ (orange dot). The blue and red curves show solutions of the full system near the canard explosion. Forcing amplitude $A$ for the blue trajectory is larger than for red and both exponentially close the canard explosion at $A^\ast$. Solutions obtained for $\eps=10^{-3}$.}
	\label{fig:spike-adding_1em3_big}
\end{figure}

Taking the fast subsystem LC family $\Gamma^{\rm a, r}$ into account, a remarkable feature of the dynamics can be seen in $(I_1, v, d)$ space (panel (d$_1$)). The spikes of the burst appear to follow the family of unstable limit cycles, thus evolving near the repelling surface $\LCmfr$. Unexpectedly, after piercing through the critical manifold, the bursting solution stays in the proximity of $\LCmfr$, instead of being repelled from it. As $I_1$ slowly drifts towards smaller values, the trajectory remains close to $\LCmfr$. These windings around $\LCmfr$ correspond to the first not fully grown spikes of the full system solution. The spikes increase in amplitude, as $I_1$ decreases, but remain small; see panel (b$_1$) and inset. An enlargement of two full system trajectories in $(I_1,v,x)$-space is shown in \cref{fig:spike-adding_1em3_big}, with two exponentially close $A$ values near $A^\ast$.

Finally, the fast oscillations terminate via an escape from $\LCmfr$. By approaching the bottom sheet of $\cmf$, the dynamics change to that of a drifting equilibrium, passing from burst to quiescence.

\subsubsection*{Explosivity and spike-adding $\boldsymbol{({\rm b}_2, {\rm c}_2,{\rm d}_2)}$.}
At a slightly larger $A$ value, close to $A^\ast$, a majority of the trajectory remains essentially unchanged, with respect to the previous $A$. The slowly drifting part along the bottom equilibrium branch of the fast subsystem and the canard segment, as well as the first oscillations, appear frozen. This is clear by comparing the blue trajectory with the red dashed curves in panels (${\rm b}_2, {\rm c}_2,{\rm d}_2$), as well as in \cref{fig:spike-adding_1em3_big}.
The fact that part of the trajectory near the fold freezes, while the following part changes significantly, is a strong indication for explosivity of the solution, when varying the parameter. This strong sensitivity towards parameter changes is typical for canard dynamics. It is caused by the presence of repelling objects in the fast subsystem, typically, but not exclusively, equilibria, like the middle sheet of $\cmf$.

Indeed, in panels (${\rm b}_2, {\rm c}_2,{\rm d}_2$), the full system trajectory possesses a canard segment staying near the middle sheet of $\cmf$, when it turns around the folded saddle $p_1$. Hence the sensitivity is to be expected. Moreover, during the fast oscillations, the bursting solution evolves close to $\LCmfr$, thus adding an additional layer of sensitivity.

Compared to the previous case in \cref{fig:spike-adding_1em3}(${\rm b}_1, {\rm c}_1,{\rm d}_1$), the full system solution winds around $\LCmfr$ more often until smaller values of $I_1$, before jumping back to $\cmf$. The trajectory essentially remains close to $\LCmfr$, but reaches up higher. This way, by passing through $A^\ast$, more and more spikes with increasing amplitude are added to the burst. These spikes are yet not fully grown to the amplitude of the stable limit cycles present in this $I_1$ region. Furthermore, we note that distance of the blue trajectory to $\LCmfr$ increases, as it winds around it, indicating some extent of repulsion near the surface; see also panel (b$_2$). This way, the full dynamics starts to escape from $\LCmfr$ and gets attracted to $\LCmfa$.

\subsubsection*{Emergence of bursting $\boldsymbol{({\rm b}_3, {\rm c}_3,{\rm d}_3)}$.} As $A$ increases further, the point at which the trajectory starts to escape from $\LCmfr$ shifts towards the lower Hopf bifurcation $\HBL$. In fact, the last two spikes are already repelled sufficiently to evolve close to the attracting surface $\LCmfa$; see panels (b$_3$,d$_3$). In other words, the number of revolutions near $\LCmfr$ reduces, while the ones around $\LCmfa$ increases. These oscillations near $\LCmfa$ are of large amplitude and mark the start of a dense burst following for larger $A$.
 
\subsubsection*{Bursting $\boldsymbol{({\rm b}_4, {\rm c}_4,{\rm d}_4)}$.} At the next step in panels ${\rm b}_4$ - ${\rm d}_4$, most of the windings around $\LCmfr$ have vanished and the fast oscillations take place in the proximity of $\LCmfa$. Additionally, the burst consists of more spikes in total, with respect to the first considered $A$ value (panels (${\rm b}_1, {\rm c}_1,{\rm d}_1$)).

Here, an additional spike-adding mechanism beyond the critical value $A^\ast$ acts and is related to the period $T(I_1)$ of $\LCmfa$ as a function of $I_1$. When $A$ is beyond the canard explosion at $A^\ast$ the full system dynamics approach $\LCmfa$ right after surpassing the lower fold $\FL$, without the excursion on $\LCmfr$. As the amplitude $A$ is increased, this attraction to $\LCmfa$ occurs at larger $I_1$ values.	The period $T(I_1$) of the fast subsystem limit cycles decreases with increasing $I_1$ and this finally leads to more windings around $\LCmfa$. On top of this, the fraction of time for which the trajectory remains near $\LCmfa$ increases. Both effects add more spikes to the burst and result in the spike-adding arches observed in \cref{fig:spike-adding_1em3}(a).

The presented results already show the complexity of how bursts are generated, that is, via a transition through the canard explosion at $A=A^\ast$, which rather surprisingly leads the canard segment to evolve around the repelling object $\LCmfr$. This is followed by spike growth via repulsion from $\LCmfr$ and attraction to $\LCmfa$, until all oscillations evolve near $\LCmfa$.

\subsection{Continuous route to bursting: spike-adding via mixed-type-like torus canards}
Before the bursting transition, the full system dynamics can be described by a single slow frequency, determined by $\eps$. After the transition, a full cycle consists of a slow phase followed by fast spiking. It is therefore characterized by the slow frequency and a fast one, given by properties of the fast subsystem.
In the context of bursting and in slow-fast systems, whose fast subsystem has both stable and unstable cycles, this change of the dynamics, from one to two frequencies, hints at a torus (Neimark-Sacker) bifurcation in the full system.
%	This change of the dynamics, from one to two frequencies, hints at a torus (Neimark-Sacker) bifurcation in the full system. In the context of bursting and in slow-fast systems whose fast subsystem has both stable and unstable cycles, then NS can play an important role

%	\mdnote[inline]{this is a bit far-fetched: spike-adding can organise the transition from rest to bursting without any NS bifurcation ... we can say that, in the context of bursting and in slow-fast systems whose fast subsystem has both stable and unstable cycles, then NS can play an important role}

This can indicate the existence of mixed type-torus canards (MTTCs). Indeed, the full system dynamics nearby the canard explosion not only exhibits a canard segment along the repelling sheet of $\cmf$, but as well a canard segment on the repelling higher dimensional invariant set $\LCmfr$. This clearly resembles mixed-type canards as described in \cite{desrochesCanardsMixedType2012}. Very similar MTTCs have been reported in \cite{baspinarCanonicalModelsTorus2021}, where they evidently arise in the singular limit $\eps \rightarrow 0$. In particular Fig. 2 of \cite{baspinarCanonicalModelsTorus2021} reports dynamics where a quasi-static motion of the full system along the attracting sheet of $\cmf$ connects to a repelling set of limit cycles created by a subcritical Hopf bifurcation.
In the NMSTP however, the understanding of mixed-type torus canards is more complex for various reasons and as we will show, only observed for small, but large enough $\eps$.

First of all, in the $\eps$ regime for which we observe mixed-type torus canards, the timescale separation is small enough for the canard segment on the middle sheet of $\cmf$ to persist as a strongly perturbed version of its singular counterpart. One can observe a turn around the lower fold $\FUc$, mediated by the folded-saddle singularity $p_1$. It forces the trajectory to pierce through $\cmf$, bringing it very close to $\LCmfr$; see \cref{fig:spike-adding_1em3_big}.

Secondly, the solutions $\Gamma^{r} \subset  \LCmfr$, despite being globally repelling, possess two stable Floquet multipliers. We argue that the associated stable directions, similar to the jump-on canard case (see \cref{sec:jump-on}), form due to the intrinsic timescales of the fast subsystem. Consequently, this means that stable manifolds can be associated with $\LCmfr$: it can attract trajectories along certain directions.

Thirdly, for large enough $\eps$, one can expect dynamics very distinct from singular orbits. In particular, solutions may not only evolve nearby, but also switch between different attracting branches of the fast subsystem. We observe this transition from the middle sheet of $\cmf$ to $\LCmfr$. Despite both being repulsive, the reasoning holds: the folded-saddle canard dynamics enforces the full system to stay close to the middle sheet of $\cmf$; stable directions allow an attraction towards $\LCmfr$; the large $\eps$ allows to bridge $\cmf$ to $\LCmfr$ and finally torus canard dynamics allow to follow $\LCmfa$ closely, adding more and more spikes in the transition region around $A^\ast$.

\subsection{Discontinuous route to bursting: block evoked by jump-on canards}\label{sec:blocking}
The transition from subthreshold oscillations to bursting when increasing the forcing amplitude $A$ explains the emergence of the very first spikes in the burst, which occur for forcing amplitudes exponentially close to $A^\ast$. They also show how the subsequent spike-adding process, reflected by the arches in \cref{fig:spike-adding_1em3}(a), occurs. In the following we will extend the analysis of this transition, taking into account different values of the parameter $\eps$. As an outlook, we will describe how, for $\eps$ values closer to the singular limit, jump-on canards interfere and block the transition.

\begin{figure}[h!]
	\includegraphics[width=1\linewidth]{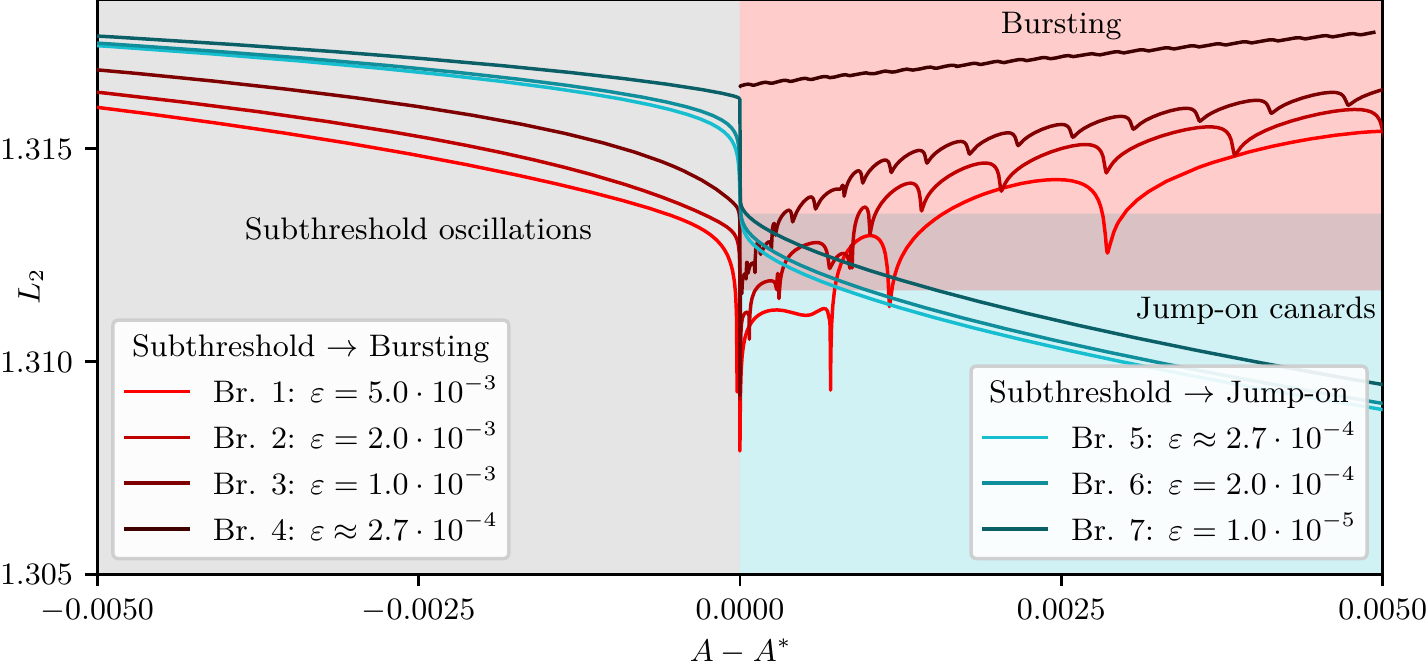}
	\caption{Families of bursting solutions and canards: The bifurcation diagram shows solution families of the full system in terms of the $L_2$-norm vs. the shifted forcing amplitude $A-A^\ast$. For all branches (br.) apart from br. 4, the $A^\ast$ denotes the location of the canard explosion. For br. 4 on the other hand it marks the termination of continuation due to insufficient accuracy. There are two types of solution families: the ones which undergo a continuous transition from subthreshold oscillations to bursting (red curves) and those which transition from subthreshold to jump-on canard dynamics (cyan curves). Br. 4 and 5 have identical $\eps\approx 2.667521298 \cdot 10^{-4}$.}
	\label{fig:diffepsibranchesl2-vs-a}
\end{figure}

In \cref{fig:diffepsibranchesl2-vs-a} solution families of periodic orbits are displayed, obtained via numerical continuation of the full system \cref{eq:model1,eq:forcing1}. The figure comprises seven branches for $\eps$ values ranging from $\eps=5\cdot10^{-3}$ to $\eps = 1\cdot10^{-5}$ and they are aligned to the canard explosion at $A=A^\ast$ (apart from branch 4). Branch 3 is identical to the one shown in \cref{fig:spike-adding_1em3} and undergoes a continuous transition from subthreshold oscillations to bursting.

As a general result, we find two types of solution families, which,  beyond the canard explosion, i.e, for $A>A^\ast$, take different paths. For $\eps \gtrapprox 2.7\cdot 10^{-4}$  the branches display a continuous transition from subthreshold oscillations to bursting; see red curves in \cref{fig:diffepsibranchesl2-vs-a}.
For $\eps \lessapprox 2.7\cdot 10^{-4}$ on the other hand the branches evolve into families of jump-on canards (curves in \cref{fig:diffepsibranchesl2-vs-a}). A value of $\eps\approx 2.667521298 \cdot 10^{-4}$ separates the two $\eps$ regimes.

However, it is clear by considering branches 4 and 5 that bursting solutions do not cease to exist for smaller $\eps$ values. Instead they coexist with the jump-on type solutions.

It remains unclear how bursting forms for $\eps \lessapprox 2.7\cdot 10^{-4}$. Nevertheless the role of jump-on canards for the emergence of bursting becomes evident: for singular or small enough $\eps \lessapprox 2.7\cdot 10^{-4}$ regular canards are observed, evolving on the middle sheet of $\cmf$. Beyond the maximal canard, the intrinsic timescales of the fast subsystem come into play and lead to the emergence of jump-on canards (see \cref{sec:jump-on}). They block a transition towards $\LCmf$ and as a consequence, bursting remains absent for these solutions families.

For $\eps \gtrapprox 2.7\cdot 10^{-4}$ however, jump-on canards cease to exist. In the amplitude regime where they would be expected, the system approaches $\LCmfa$ and bursts instead. The region around the canard explosion is populated by MTTCs and separates subthreshold oscillations from the bursting regime.
It remains an open question for future work how the differentiation depending on $\eps$ occurs and in particular which possible bifurcations of the full dynamics result in the distinct regimes of continuous bursting transition and blocking jump-on canards.
%\mdnote[inline]{we'll add something about future work here...}

\section{Network behaviour}\label{sec:netbehaviour}
The neural mass model with STP \cref{eq:model1} is an exact limit of the underlying QIF network \cref{eq:QIFnetwork2} as $N\rightarrow \infty$. We want to emphasize the benefits of neural mass models and their capability of describing neuronal dynamics at a macroscopic scale. 
In \cref{fig:bifdiag}(b-f) results using the fast subsystem and the corresponding QIF network have been shown. For the original neural mass model \cite{montbrioMacroscopicDescriptionNetworks2015} the exactness of the meanfield limit has been exploited in various studies in order to understand the collective dynamics of large neuronal populations \cite{pazoQuasiperiodicPartialSynchronization2016, devalleFiringRateEquations2017a,	schmidtNetworkMechanismsUnderlying2018b,	ceniCrossFrequencyCoupling2020, segneriThetanestedGammaOscillations2020a, taherExactNeuralMass2020}. We want to note that STP, as opposed to exponential synapses used in Ref. \cite{avitabileCrossscaleExcitabilityNetworks2021}, results in a substantially higher sensitivity towards finite size fluctuations and numerical errors, rendering the agreement of network and neural mass less clear, in particular in the canard regime.

Here we want to assess if the mechanisms leading to bursting in the neural mass model persists in a finite-sized network. To our
knowledge such analysis, in particular in presence of short-term synaptic plasticity, has not been performed.
For this we introduce the external periodic forcing $I_1=A\sin (\eps t)$ also into the network  \cref{eq:QIFnetwork2} and investigate the QIF population dynamics nearby the canard explosion of the neural mass.
\begin{figure}[h!]
	\includegraphics[width=1\linewidth]{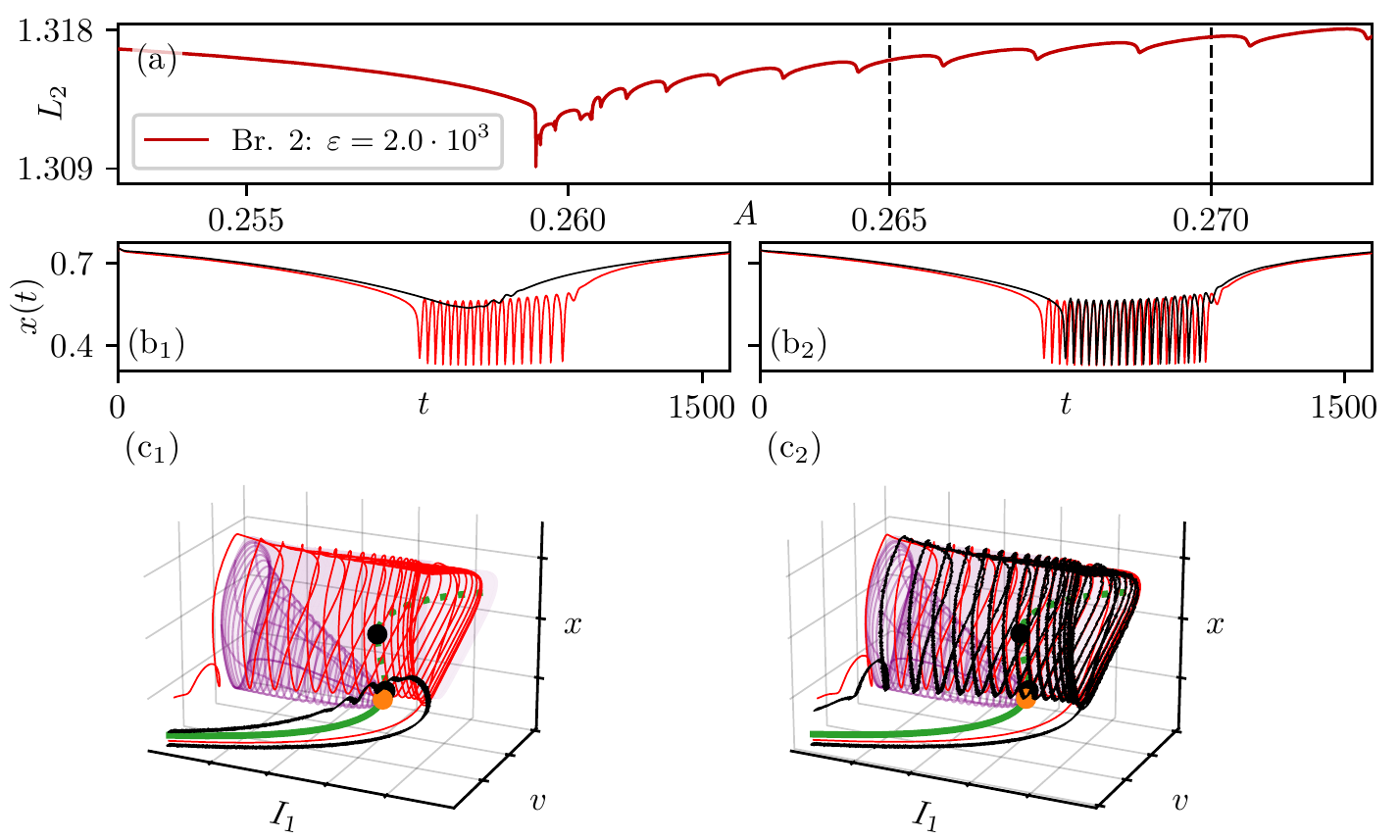}
	\caption{Bursting in the QIF network: {\bf (a)} The bifurcation diagram shows a family of periodic orbits of the full neural mass system in terms of the $L_2$-norm vs. the forcing amplitude $A-A^\ast$. It corresponds to br. 2 in \cref{fig:diffepsibranchesl2-vs-a}. The black dashed lines mark $A=0.265$ and $A=0.270$ respectively. {\bf (b - c)} Trajectories of the full system (red) superimposed on results obtained by simulating the QIF network \cref{eq:QIFnetwork2} with $N=100000$ neurons and in presence of sinusoidal forcing $I_1(t)=A\sin(\eps t)$. In (b$_1$,c$_1$) the forcing amplitude is given by $A=0.265$;  in (b$_2$,c$_2$) by $A=0.270$. Row (b) shows the time series of the solution in terms of $x(t)$ vs. time $t$. In row (c) they are shown in $(I_1,v,x)$-space together with the critical manifold $\cmf$, with attracting (repelling) sheets as a solid (dashed) green line, and the invariant manifold $\LCmf$ of the fast subsystem (purple wireframe and surface). The orange dot denotes $\HBLc$, the black dots $\FLc$ and $\FUc$. The $x$-axis is inverted in (c).
	}
	\label{fig:burst_network_2em3}
\end{figure}
In \cref{fig:burst_network_2em3}(a) the family of periodic solutions transitioning from subthreshold oscillation to bursting is shown for $\eps=2\cdot 10^{-3}$. Two dashed black lines mark the values $A=0.265$ and $A=0.27$, respectively, for which neural mass  (red) as well as QIF network trajectories (black) are depicted in \cref{fig:burst_network_2em3}(b,c). Row (b) shows the time series $x(t)$ vs. $t$, row (c) the trajectories in $(I_1,v,x)$-space.

The QIF network consists of $N=100000$ neurons and $\eps=2\cdot 10^{-3}$ is chosen to maintain reasonable computation times.
In the network, the initial conditions are chosen according to fixed point values $(r^\ast, v^\ast,x^\ast,u^\ast)$ obtained from the neural mass at $I_1=0$. However, initializing the network firing rate and average membrane potential at a given value is a non-trivial problem. To overcome this, the initial conditions of the QIF network are set to $V_i=v^\ast\quad \text{for all }i=1,\dots,N$, $x=x^\ast$ and $u=u^\ast$. After a short transient, at time $t=0$, the network has reached equilibrium and the forcing $I=A\sin (\eps t)$ sets in.

We start the analysis by considering $A=0.265$, marked by the left most dashed line in \cref{fig:burst_network_2em3}(a). It lays within the spike-adding arches of the neural mass, which displays bursting with a multitude of fully grown spikes; see red curves in \cref{fig:burst_network_2em3}(b$_1$, c$_1$). They mainly evolve in the proximity of the invariant set $\LCmfa$. This periodic solution of the full problem does not highlight a canard segment around the folded singularity $p_1$, which is to be expected, since $A$ is not close enough to the canard explosion at $A^\ast$.

On the other hand, the macroscopic state of the QIF network is found to be just at the start of the spike-adding process. The orbit possesses a canard segment: the turn around the lower fold $\FL$ is evident. This is clear by looking at the black curve in \cref{fig:burst_network_2em3}(c$_1$). More strikingly, a few low-amplitude oscillations are picked up in the vicinity of the repelling invariant set $\LCmfr$. This means that, despite the discrepancies, the mixed-type torus canards observed in the neural mass are also found in the QIF network. It is plausible to assume that spike-adding in the network follows the same mechanism as described in \cref{sec:transbursting}.

At the larger amplitude $A=0.270$, shown in \cref{fig:burst_network_2em3}(b$\blue{_2}$, c$\blue{_2}$), the agreement between neural mass and network  trajectory is more distinct. In both systems the canard segment is absent and they both exhibit bursting with a comparable number of large amplitude spikes. Taking also into account the previous case $A=0.265$, it appears that the bifurcation structure of the network is shifted towards larger amplitudes with respect to the neural mass.

\section{Discussion}

The results of this work lay within the intersection of various fields including meanfield theory and slow-fast dynamics.
First of all, we make use of recent developments of meanfield theory, namely the powerful OA ansatz, in order to understand the emergent dynamics of spiking neuronal networks on macroscopic scale.

Secondly, we aim at more biologically plausible models in this work, through the inclusion of synaptic dynamics, in form of STP, while preserving the exactness of the meanfield limit with respect to the QIF network. This adds more realism to the QIF network and neural mass, but at the same time adds to the complexity of the collective dynamics, even in absence of external forcing. In the part of parameter space chosen here, STP allows for the existence limit cycles, that are absent in the original model. As we have shown via extensive numerical evidence, STP leads to various peculiarities, when forcing the system slowly and periodically.

That leads to the third point, which we put our focus on: the treatment of the forced neural mass with STP using methods of singular perturbation theory, in particular slow-fast dissection. Without STP, slow external forcing already gives rise to bursts, as shown in \cite{montbrioMacroscopicDescriptionNetworks2015}. The fast oscillations of these orbits vanish in the singular limit and relaxation oscillations remain. Accounting for STP, however, leads to more intricate dynamics.

One of the fundamental elements are canards. Expectedly, due to the slow harmonic passage through a fold of the S-shaped critical manifold, they appear as folded-saddle canards, which play a role for spike-adding in parabolic bursters. Here, despite the fact that the equilibrium destabilizing bifurcation is a subcritical Hopf bifurcation, the observed bursts are reminiscent of elliptic bursting (i.e., subcritical fold/fold of cycles bursting), due the slow passage through the Hopf bifurcation.

Concerning the full system dynamics, a folded saddle can be found, with associated explosive canard solutions. It separates quiescent orbits with purely slow dynamics from bursting ones. In this transition region we find an intriguing interplay of slow-fast effects: jump-on canards exist close enough to the singular limit and are associated with a subtle timescale separation of the neural mass, allegedly invoked by STP. They connect two repelling sheets of the critical manifold and more strikingly, block a continuous transition from quiescence to bursting.

Jump-on canards are one of the surprising elements of this work. However, they vanish when considering biologically more plausible frequencies of the periodic external current. Despite the rather insufficient timescale separation in this parameter regime, we find orbits which display a strongly perturbed canard segment. It is remarkable that these canards, as opposed to the jump-on case, do not block a continuous transition towards bursting, but much more promote it. Specifically they guide the trajectories into the proximity of unstable limit cycles of the fast subsystem. Here once again, an unexpected mechanisms sets in and allows attraction towards these globally repelling cycles. The final element to the spike-adding mechanisms in this region are rapidly oscillating segments that revolve around the branch of unstable limit cycles and consecutively add more spikes to the orbit.

Overall, the full dynamics nearby the canard explosion can be seen as mixed-type-like canards of peculiar nature: orbits evolve nearby the repelling middle sheet of the critical manifolds, as well nearby unstable limit cycles. Typically mixed-type canards are observed towards the singular limit, as they represent a slow-fast effect. The neural mass with STP exhibits mixed-type-like canards only for large enough forcing frequency; the phenomenon disappears for too slow forcing and is blocked by jump-on canards.

By virtue of the NMSTP model being the main subject of this work, we want to emphasize that it represents an exact limit of the QIF network with STP. Discrepancies between network and NMSTP arise due to numerical errors and finite-size fluctuations, especially under slow forcing. Nevertheless, our results clearly show a good agreement between the two models. In particular, the network simulations display the mixed-type-like torus canard dynamics, hence this can be regarded as a strong evidence for the same mechanisms to be responsible for burst spike-adding in the network.

In summary, the NMSTP turns out to be a useful approach in order to investigate the ensemble dynamics of neuronal populations in presence of STP. Slow-fast dissection reveals the mechanisms underlying burst generation on population level. As is it turns out, synaptic dynamics indeed enriches the complexity of the problem, by giving rise to peculiar jump-on and mixed-type like torus canards, which both appear due to the timescales associated with STP.

Finally, while this work is in the scope of neural mass models and STP, the methodology of OA reduction and slow-fast dissection, coupled with numerical bifurcation analysis, can be applied to a much broader class of phase oscillator systems, like the Kuramoto model. This can lead to a better understanding of emerging collective slow-fast dynamics of large scale networks.

\cleardoublepage

\color{black}
\bibliography{ref_final.bib}

\end{document}